\documentclass{article}
\usepackage{amssymb,amsmath,theorem,euscript}

\input{epsf}

\newcounter{sec}

\def\sm{\smallskip}

\newcounter{punct}[sec]

\def\punct{\refstepcounter{punct}{\arabic{sec}.\arabic{punct}.  }}

\def\COUNTERS{\addtocounter{sec}{1}
              \setcounter{punct}{0}
          \setcounter{equation}{0}
          \setcounter{theorem}{0}
                  }
\newtheorem{theorem}{Theorem}[sec]
\newtheorem{proposition}[theorem]{Proposition}
\newtheorem{lemma}[theorem]{Lemma}
\newtheorem{corollary}[theorem]{Corollary}

\newtheorem{definition}[theorem]{Definition}

\begin{document}

 \def\ov{\overline}
\def\wt{\widetilde}
 \newcommand{\rk}{\mathop {\mathrm {rk}}\nolimits}
\newcommand{\Aut}{\mathop {\mathrm {Aut}}\nolimits}
\newcommand{\Out}{\mathop {\mathrm {Out}}\nolimits}
\renewcommand{\Re}{\mathop {\mathrm {Re}}\nolimits}
\renewcommand{\Im}{\mathop {\mathrm {Im}}\nolimits}
 \newcommand{\tr}{\mathop {\mathrm {tr}}\nolimits}
  \newcommand{\Hom}{\mathop {\mathrm {Hom}}\nolimits}
   \newcommand{\diag}{\mathop {\mathrm {diag}}\nolimits}
   \newcommand{\supp}{\mathop {\mathrm {supp}}\nolimits}
 \newcommand{\im}{\mathop {\mathrm {im}}\nolimits}

\def\Br{\mathrm {Br}}

\def\SL{\mathrm {SL}}
\def\SU{\mathrm {SU}}
\def\GL{\mathrm {GL}}
\def\U{\mathrm U}
\def\OO{\mathrm O}
 \def\Sp{\mathrm {Sp}}
 \def\SO{\mathrm {SO}}
\def\SOS{\mathrm {SO}^*}
 \def\Diff{\mathrm{Diff}}
 \def\Vect{\mathfrak{Vect}}
\def\PGL{\mathrm {PGL}}
\def\PU{\mathrm {PU}}
\def\PSL{\mathrm {PSL}}
\def\Symp{\mathrm{Symp}}
\def\End{\mathrm{End}}
\def\Mor{\mathrm{Mor}}
\def\Aut{\mathrm{Aut}}
 \def\PB{\mathrm{PB}}
 \def\cA{\mathcal A}
\def\cB{\mathcal B}
\def\cC{\mathcal C}
\def\cD{\mathcal D}
\def\cE{\mathcal E}
\def\cF{\mathcal F}
\def\cG{\mathcal G}
\def\cH{\mathcal H}
\def\cJ{\mathcal J}
\def\cI{\mathcal I}
\def\cK{\mathcal K}
 \def\cL{\mathcal L}
\def\cM{\mathcal M}
\def\cN{\mathcal N}
 \def\cO{\mathcal O}
\def\cP{\mathcal P}
\def\cQ{\mathcal Q}
\def\cR{\mathcal R}
\def\cS{\mathcal S}
\def\cT{\mathcal T}
\def\cU{\mathcal U}
\def\cV{\mathcal V}
 \def\cW{\mathcal W}
\def\cX{\mathcal X}
 \def\cY{\mathcal Y}
 \def\cZ{\mathcal Z}
\def\0{{\ov 0}}
 \def\1{{\ov 1}}
 \def\frA{\mathfrak A}
 \def\frB{\mathfrak B}
\def\frC{\mathfrak C}
\def\frD{\mathfrak D}
\def\frE{\mathfrak E}
\def\frF{\mathfrak F}
\def\frG{\mathfrak G}
\def\frH{\mathfrak H}
\def\frI{\mathfrak I}
 \def\frJ{\mathfrak J}
 \def\frK{\mathfrak K}
 \def\frL{\mathfrak L}
\def\frM{\mathfrak M}
 \def\frN{\mathfrak N} \def\frO{\mathfrak O} \def\frP{\mathfrak P} \def\frQ{\mathfrak Q} \def\frR{\mathfrak R}
 \def\frS{\mathfrak S} \def\frT{\mathfrak T} \def\frU{\mathfrak U} \def\frV{\mathfrak V} \def\frW{\mathfrak W}
 \def\frX{\mathfrak X} \def\frY{\mathfrak Y} \def\frZ{\mathfrak Z} \def\fra{\mathfrak a} \def\frb{\mathfrak b}
 \def\frc{\mathfrak c} \def\frd{\mathfrak d} \def\fre{\mathfrak e} \def\frf{\mathfrak f} \def\frg{\mathfrak g}
 \def\frh{\mathfrak h} \def\fri{\mathfrak i} \def\frj{\mathfrak j} \def\frk{\mathfrak k} \def\frl{\mathfrak l}
 \def\frm{\mathfrak m} \def\frn{\mathfrak n} \def\fro{\mathfrak o} \def\frp{\mathfrak p} \def\frq{\mathfrak q}
 \def\frr{\mathfrak r} \def\frs{\mathfrak s} \def\frt{\mathfrak t} \def\fru{\mathfrak u} \def\frv{\mathfrak v}
 \def\frw{\mathfrak w} \def\frx{\mathfrak x} \def\fry{\mathfrak y} \def\frz{\mathfrak z} \def\frsp{\mathfrak{sp}}
 \def\bfa{\mathbf a} \def\bfb{\mathbf b} \def\bfc{\mathbf c} \def\bfd{\mathbf d} \def\bfe{\mathbf e} \def\bff{\mathbf f}
 \def\bfg{\mathbf g} \def\bfh{\mathbf h} \def\bfi{\mathbf i} \def\bfj{\mathbf j} \def\bfk{\mathbf k} \def\bfl{\mathbf l}
 \def\bfm{\mathbf m} \def\bfn{\mathbf n} \def\bfo{\mathbf o} \def\bfp{\mathbf p} \def\bfq{\mathbf q} \def\bfr{\mathbf r}
 \def\bfs{\mathbf s} \def\bft{\mathbf t} \def\bfu{\mathbf u} \def\bfv{\mathbf v} \def\bfw{\mathbf w} \def\bfx{\mathbf x}
 \def\bfy{\mathbf y} \def\bfz{\mathbf z} \def\bfA{\mathbf A} \def\bfB{\mathbf B} \def\bfC{\mathbf C} \def\bfD{\mathbf D}
 \def\bfE{\mathbf E} \def\bfF{\mathbf F} \def\bfG{\mathbf G} \def\bfH{\mathbf H} \def\bfI{\mathbf I} \def\bfJ{\mathbf J}
 \def\bfK{\mathbf K} \def\bfL{\mathbf L} \def\bfM{\mathbf M} \def\bfN{\mathbf N} \def\bfO{\mathbf O} \def\bfP{\mathbf P}
 \def\bfQ{\mathbf Q} \def\bfR{\mathbf R} \def\bfS{\mathbf S} \def\bfT{\mathbf T} \def\bfU{\mathbf U} \def\bfV{\mathbf V}
 \def\bfW{\mathbf W} \def\bfX{\mathbf X} \def\bfY{\mathbf Y} \def\bfZ{\mathbf Z} \def\bfw{\mathbf w}
 \def\R {{\mathbb R }} \def\C {{\mathbb C }} \def\Z{{\mathbb Z}} \def\H{{\mathbb H}} \def\K{{\mathbb K}}
 \def\N{{\mathbb N}} \def\Q{{\mathbb Q}} \def\A{{\mathbb A}} \def\T{\mathbb T} \def\P{\mathbb P} \def\G{\mathbb G}
 \def\bbA{\mathbb A} \def\bbB{\mathbb B} \def\bbD{\mathbb D} \def\bbE{\mathbb E} \def\bbF{\mathbb F} \def\bbG{\mathbb G}
 \def\bbI{\mathbb I} \def\bbJ{\mathbb J} \def\bbL{\mathbb L} \def\bbM{\mathbb M} \def\bbN{\mathbb N} \def\bbO{\mathbb O}
 \def\bbP{\mathbb P} \def\bbQ{\mathbb Q} \def\bbS{\mathbb S} \def\bbT{\mathbb T} \def\bbU{\mathbb U} \def\bbV{\mathbb V}
 \def\bbW{\mathbb W} \def\bbX{\mathbb X} \def\bbY{\mathbb Y} \def\kappa{\varkappa} \def\epsilon{\varepsilon}
 \def\phi{\varphi} \def\le{\leqslant} \def\ge{\geqslant}

\def\UU{\bbU}
\def\Mat{\mathrm{Mat}}
\def\tto{\rightrightarrows}

\def\Gms{\mathrm {Gms}}
\def\Ams{\mathrm {Ams}}
\def\Isom{\mathrm {Isom}}

\def\Gr{\mathrm{Gr}}

\def\graph{\mathrm{graph}}

\def\O{\mathrm{O}}

\def\la{\langle}
\def\ra{\rangle}


 \def\ov{\overline}
\def\wt{\widetilde}

\renewcommand{\Re}{\mathop {\mathrm {Re}}\nolimits}
\def\Br{\mathrm {Br}}

 \def\Isom{\mathrm {Isom}}
\def\SL{\mathrm {SL}}
\def\SU{\mathrm {SU}}
\def\GL{\mathrm {GL}}
\def\U{\mathrm U}
\def\OO{\mathrm O}
 \def\Sp{\mathrm {Sp}}
 \def\SO{\mathrm {SO}}
\def\SOS{\mathrm {SO}^*}
 \def\Diff{\mathrm{Diff}}
 \def\Vect{\mathfrak{Vect}}
\def\PGL{\mathrm {PGL}}
\def\PU{\mathrm {PU}}
\def\PSL{\mathrm {PSL}}
\def\Symp{\mathrm{Symp}}
\def\End{\mathrm{End}}
\def\Mor{\mathrm{Mor}}
\def\Aut{\mathrm{Aut}}
 \def\PB{\mathrm{PB}}
 \def\cA{\mathcal A}
\def\cB{\mathcal B}
\def\cC{\mathcal C}
\def\cD{\mathcal D}
\def\cE{\mathcal E}
\def\cF{\mathcal F}
\def\cG{\mathcal G}
\def\cH{\mathcal H}
\def\cJ{\mathcal J}
\def\cI{\mathcal I}
\def\cK{\mathcal K}
 \def\cL{\mathcal L}
\def\cM{\mathcal M}
\def\cN{\mathcal N}
 \def\cO{\mathcal O}
\def\cP{\mathcal P}
\def\cQ{\mathcal Q}
\def\cR{\mathcal R}
\def\cS{\mathcal S}
\def\cT{\mathcal T}
\def\cU{\mathcal U}
\def\cV{\mathcal V}
 \def\cW{\mathcal W}
\def\cX{\mathcal X}
 \def\cY{\mathcal Y}
 \def\cZ{\mathcal Z}
\def\0{{\ov 0}}
 \def\1{{\ov 1}}
 \def\frA{\mathfrak A}
 \def\frB{\mathfrak B}
\def\frC{\mathfrak C}
\def\frD{\mathfrak D}
\def\frE{\mathfrak E}
\def\frF{\mathfrak F}
\def\frG{\mathfrak G}
\def\frH{\mathfrak H}
\def\frI{\mathfrak I}
 \def\frJ{\mathfrak J}
 \def\frK{\mathfrak K}
 \def\frL{\mathfrak L}
\def\frM{\mathfrak M}
 \def\frN{\mathfrak N} \def\frO{\mathfrak O} \def\frP{\mathfrak P} \def\frQ{\mathfrak Q} \def\frR{\mathfrak R}
 \def\frS{\mathfrak S} \def\frT{\mathfrak T} \def\frU{\mathfrak U} \def\frV{\mathfrak V} \def\frW{\mathfrak W}
 \def\frX{\mathfrak X} \def\frY{\mathfrak Y} \def\frZ{\mathfrak Z} \def\fra{\mathfrak a} \def\frb{\mathfrak b}
 \def\frc{\mathfrak c} \def\frd{\mathfrak d} \def\fre{\mathfrak e} \def\frf{\mathfrak f} \def\frg{\mathfrak g}
 \def\frh{\mathfrak h} \def\fri{\mathfrak i} \def\frj{\mathfrak j} \def\frk{\mathfrak k} \def\frl{\mathfrak l}
 \def\frm{\mathfrak m} \def\frn{\mathfrak n} \def\fro{\mathfrak o} \def\frp{\mathfrak p} \def\frq{\mathfrak q}
 \def\frr{\mathfrak r} \def\frs{\mathfrak s} \def\frt{\mathfrak t} \def\fru{\mathfrak u} \def\frv{\mathfrak v}
 \def\frw{\mathfrak w} \def\frx{\mathfrak x} \def\fry{\mathfrak y} \def\frz{\mathfrak z} \def\frsp{\mathfrak{sp}}
 \def\bfa{\mathbf a} \def\bfb{\mathbf b} \def\bfc{\mathbf c} \def\bfd{\mathbf d} \def\bfe{\mathbf e} \def\bff{\mathbf f}
 \def\bfg{\mathbf g} \def\bfh{\mathbf h} \def\bfi{\mathbf i} \def\bfj{\mathbf j} \def\bfk{\mathbf k} \def\bfl{\mathbf l}
 \def\bfm{\mathbf m} \def\bfn{\mathbf n} \def\bfo{\mathbf o} \def\bfp{\mathbf p} \def\bfq{\mathbf q} \def\bfr{\mathbf r}
 \def\bfs{\mathbf s} \def\bft{\mathbf t} \def\bfu{\mathbf u} \def\bfv{\mathbf v} \def\bfw{\mathbf w} \def\bfx{\mathbf x}
 \def\bfy{\mathbf y} \def\bfz{\mathbf z} \def\bfA{\mathbf A} \def\bfB{\mathbf B} \def\bfC{\mathbf C} \def\bfD{\mathbf D}
 \def\bfE{\mathbf E} \def\bfF{\mathbf F} \def\bfG{\mathbf G} \def\bfH{\mathbf H} \def\bfI{\mathbf I} \def\bfJ{\mathbf J}
 \def\bfK{\mathbf K} \def\bfL{\mathbf L} \def\bfM{\mathbf M} \def\bfN{\mathbf N} \def\bfO{\mathbf O} \def\bfP{\mathbf P}
 \def\bfQ{\mathbf Q} \def\bfR{\mathbf R} \def\bfS{\mathbf S} \def\bfT{\mathbf T} \def\bfU{\mathbf U} \def\bfV{\mathbf V}
 \def\bfW{\mathbf W} \def\bfX{\mathbf X} \def\bfY{\mathbf Y} \def\bfZ{\mathbf Z} \def\bfw{\mathbf w}
 \def\R {{\mathbb R }} \def\C {{\mathbb C }} \def\Z{{\mathbb Z}} \def\H{{\mathbb H}} \def\K{{\mathbb K}}
 \def\N{{\mathbb N}} \def\Q{{\mathbb Q}} \def\A{{\mathbb A}} \def\T{\mathbb T} \def\P{\mathbb P} \def\G{\mathbb G}
 \def\bbA{\mathbb A} \def\bbB{\mathbb B} \def\bbD{\mathbb D} \def\bbE{\mathbb E} \def\bbF{\mathbb F} \def\bbG{\mathbb G}
 \def\bbI{\mathbb I} \def\bbJ{\mathbb J} \def\bbL{\mathbb L} \def\bbM{\mathbb M} \def\bbN{\mathbb N} \def\bbO{\mathbb O}
 \def\bbP{\mathbb P} \def\bbQ{\mathbb Q} \def\bbS{\mathbb S} \def\bbT{\mathbb T} \def\bbU{\mathbb U} \def\bbV{\mathbb V}
 \def\bbW{\mathbb W} \def\bbX{\mathbb X} \def\bbY{\mathbb Y} \def\kappa{\varkappa} \def\epsilon{\varepsilon}
 \def\phi{\varphi} \def\le{\leqslant} \def\ge{\geqslant}

\def\UU{\bbU}
\def\Mat{\mathrm{Mat}}
\def\tto{\rightrightarrows}

\def\Gr{\mathrm{Gr}}

\def\graph{\mathrm{graph}}

\def\O{\mathrm{O}}

\def\la{\langle}
\def\ra{\rangle}

\def\wre{\divideontimes}

\begin{center}
\bf\Large
Infinite symmetric groups and combinatorial constructions of topological field theory type
 \\
 \medskip
 \sc\large
 Yu. A. Neretin%
 \footnote{Supported by the grant FWF, P25142.}
\end{center}

\begin{flushright}
Dedicated to A.M.Vershik
\end{flushright}

{\small The paper contains a survey of train constructions
for  infinite symmetric groups and related groups.
For certain pairs  (a group $G$, a subgroup  $K$), 
we construct categories, whose morphisms are
two-dimensional surfaces tiled by polygons and colored in a certain way.
A product of morphisms is a gluing of combinatorial
bordisms. For a unitary representation of
$G$ we assign a functor from the category of bordisms to the category of Hilbert spaces and bounded operators.
The construction has numerous variations, instead of surfaces there arise also 
one-dimensional objects of Brauer diagram type, multi-dimensional
pseudomanifolds, bipartite graphs.}

\medskip
{\sc 
\noindent
1. Introduction. The Thoma theorem and bisymmetric group.
\newline
2. The Lieberman theorem. 
\newline
3. The trisymmetric group and triangulated two-dimensional
\newline
\phantom{7.a}  bordisms.
\newline
4. One-dimensional constructions (chips). 
\newline
\phantom{7.a} Addendum. The classification of representations of 
\newline
\phantom{7.a} the  bisymmetric group.
\newline
5. Two-dimensional constructions. 
\newline
6. Categories of bipartite graphs.
\newline
7. Bordisms of pseudomanifolds.
\newline
8. Spherical functions with respect to the Young subgroup. 
\newline
\phantom{7.} The Nessonov theorem.
}

\section{Introduction}

\COUNTERS

{\bf\punct  Infinite symmetric groups and topological field
theories.} Representation theory of infinite symmetric groups  $S_\infty$ was initiated by
two 'orthogonal' works. The first one was the paper of Elmar Thoma
\cite{Tho1},
1964, where he introduced analogs of characters for the group 
 $S_\infty$ of finitely supported infinite permutations.
The second was the paper of Arthur Lieberman 
 \cite{Lie}, 1972, where he classified 
 all unitary representations of
 the complete infinite symmetric group.

 Initially, Thoma introduced characters of the group
   $S_\infty$ as traces of operators of representations in type $II_1$
   Murray--von Neumann factors. Below we mention the term 'factors' several times,
   but actually these objects are not used below and are not necessary for understanding of the paper.
   Further, A.M.Vershik and S.V.Kerov 
    \cite{VK1} constructed explicitly the representations of $S_\infty$ 
 in factors corresponding to the Thoma characters (see, also \cite{Ver}),
 and showed how to obtain Thoma characters as limits of characters of finite symmetric groups
 \cite{VK2}. 
 
 In 1989, G.I.Olshanski
 \cite{Olsh-symm} noticed that the Thoma characters correspond 
 to representations of the double  $S_\infty\times S_\infty$ having fixed vectors
 with respect to the diagonal subgroup. In the same work, there was obtained
  a first example of (one-dimensional) topological  field theory from
  representations of infinite symmetric groups (sense of this sentence is explained below).
One of recent achievements of representation theory of infinite symmetric groups
is the explicit decomposition of an analog of regular representations for
the double
$S_\infty\times S_\infty$  in \cite{KOV}.

Next, let us say about the term 'topological field theory' in the sense of the definition
proposed by M.Atiyah in
\cite{Ati} (see also a recent survey  \cite{Tel}).
Consider a category, whose objects are 
 $(k-1)$-dimensional manifolds, a morphism $M$ to $N$ is a  $k$-dimensional
 manifold  $R$, whose boundary is identified with 
 the disjoint union 
 $M\coprod N$. The product of two morphisms 
$R:M\to N$, $Q:N\to K$ is the result of gluing  of $R$ and $Q$
along the manifold $N$. A topological field theory
is a functor from this category to a category
of linear spaces and linear operators%
\footnote{
The term 'bordisms' is used in this context.
However, the theory of bordisms is an elder domain of mathematics
 (see. e.g., \cite{Swi}), conflicts with the well-established terminology easily
 arise here. Apparently, the term 'field theory'
 was introduced by analogy 
 with 'conformal field theories', see \cite{Seg}, \cite{Ner-tube}.}%
$^,$%
\footnote{Two-dimensional 'topological field theories', where surfaces 
are tiled by polygons were considered in \cite{Bae}, \cite{Nat}.}.
The content of the present paper is constructing
of combinatorial objects of 'topological field theory' type
from representations of infinite symmetric groups
(\cite{Olsh-symm}, \cite{Ner-symm}, \cite{Ner-preprint},  \cite{GN}).

This paper is more-or-less self-closed, a minimal familiarity of the reader with representation
theory is assumed.

\sm

{\bf\punct The infinite symmetric group.%
\label{ss:group}} By $S_\infty$ we denote the group of all finitely supported
permutations of a countable set.
On default we assume that the countable set is the set $\N$.
If we need for another countable set 
 $\Omega$, we write   $S_\infty(\Omega)$. 
The group  $S_\infty$ is the inductive limit of finite symmetric groups
 $S_n$,
$$
S_\infty:=\lim_{\longrightarrow} S_n.
$$
The group $S_\infty$ can be realized as the group
of infinite invertible  0-1-matrices%
\footnote{By a  {\it 0-1-matrix} we call a matrix such that all matrix elements
are  0 or 1 and each column and each row contains at most one unit.
 Such matrix is invertible iff each column and each row contains 
precisely one unit.}
$g$ such that $g_{kk}=1$ for sufficiently large $k$.
Also, we will represent permutations (finite or infinite) by pictures 
of the type
$$
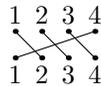
\epsfbox{bord.2}
$$

The group $S_\infty$ is a countable discrete group,
it is not a type $I$  group
(see, e.g., \cite{Dix}, \cite{Tho2}),
and a representation theory in the usual sense for this group is impossible.

Let $\alpha=0$, $1$, $2$, \dots
Denote by $S_\infty[\alpha]$ the subgroup in  $S_\infty$, consisting of permutations
fixing the points
$1$, \dots, $\alpha$. This group consists of  0-1 matrices of the form 
\begin{equation}
g=\begin{pmatrix}
                                                                 1_\alpha&0\\0&r
                                                                \end{pmatrix}, 
                                                                \text{where $r\in S_\infty$}
                                                                \label{eq:Salpha}
                                                                ,\end{equation}
    and $1_\alpha$ is the unit matrix of size $\alpha$.  
    We get a sequence of subgroups
$$
S_\infty=S_\infty[0]\supset S_\infty[1]\supset S_\infty[2]\supset\dots
$$ 
All subgroups 
$S_\infty[\alpha]$ are canonically isomorphic
to the group $S_\infty$, in notation (\ref{eq:Salpha}) the isomorphism
is done by
 $r\mapsto g$.

\sm

{\bf\punct Complete infinite symmetric group.%
\label{ss:complete-group}} By $\ov S_\infty$ we denote the group of all permutations 
of a countable set. By $\ov S_\infty(\Omega)$ denote the group of all permutations
of a countable set  $\Omega$. As above, we define the subgroups   $\ov S_\infty[\alpha]$
consisting of matrices of the form
$ \begin{pmatrix} 1_\alpha&0\\0&r \end{pmatrix}$, where $r\in \ov S_\infty$.

We define a totally disconnected topology on $S_\infty$ 
assuming that the subgroups  $\ov S_\infty[\alpha]$ are open.
We give two equivalent definitions of this topology:

-- $g_j\to g$ in $\ov S_\infty$ if for any  $k\in\N$ starting some $j$ we have 
$g_j (k)=g(k)$;

-- $g_j\to g$ if the sequence of corresponding 0-1-matrices weakly converges%
\footnote{See, e.g., \cite{KG}. In our case, this means 
the element-wise convergence of matrices, i.e.,
a stabilization.}.

The group $\ov S_\infty$ is complete  (in the Raikov sense%
\footnote{I.e., any two-side-fundamental sequence converges.
A sequence  $g_j$ is called two-side-fundamental if the both double sequences
 $g_j g_k^{-1}$ and $g_k^{-1}g_j $ tend to unit if $j$, $k\to\infty$.
 This definition differs from Bourbaki's definition of completeness. Bourbaki requires
 the convergence of left fundamental sequences and right-fundamental sequences, see
 \cite{Rai}, \cite{Bou}.}) topological group with respect 
 to our topology%
\footnote{%
This is a unique separable topology on   $\ov S_\infty$ (A.Rosendal).
Recall that a topological group is called Polish, if it is metrizable,
separable, and complete
(another version of the definition: if the topology is Polish). For Polish groups
there is  a wide collection of theorems about automatic continuity of measurable homomorphisms,
about automatically continuity of arbitrary homomorphisms, about uniqueness of Polish topologies,
see
 \cite{Kech}, \cite{Tsa} and references in the last work.}.

In 1972 A.Lieberman
 \cite{Lie} obtained a complete classification 
 of unitary representations of the group   $\ov S_\infty$.
 This classification is uncomplicated: all irreducible
 representations are induced from representations  $\rho\otimes I$
of subgroups $S(\alpha)\times \ov S_\infty[\alpha]$, where $\rho$ is irreducible and 
$I$ is the trivial one-dimensional representation of 
$\ov S_\infty[\alpha]$ (for further details, see below \S2).     
On the other hand this classification, in a certain sense, 
is not too interesting, since the representation theory $\ov S_\infty$
contains no new essences with respect to the representation theory of
finite symmetric groups.

However, a representation theory of infinite symmetric group exists, but it 
is neither a representation theory of 
$S_\infty$, nor a representation theory of $\ov S_\infty$. A key step
was the work of Thoma  \cite{Tho1}, 1964.

\sm

{\bf\punct The Thoma theorem.%
\label{ss:thoma}} Let $f$ be a complex-valued function on a countable
group  $G$. It is called {\it central} if it is constant on conjugacy classes. 
It is {\it positive definite} if for any finite collection 
$g_j\in G$, the matrix composed of    $f(g_ig_j^{-1})$ is positive semi-definite.
The set of central positive definite functions is a convex cone. 
By $\cK(G)$ we denote the set of such functions satisfying the 
condition  $f(1)=1$. 
The set $\cK(G)$ is convex and compact with respect to 
the topology of point-wise convergence. By the Krein--Milman theorem
 (see, e.g.,  \cite{KG}), $\cK(G)$ is the convex hull of its set of extreme points.

\begin{definition}
{\rm(\cite{Tho1})} A character of the group $G$ is an extreme point of the set $\cK(G)$.
\end{definition}

Let us remark, that for finite or compact groups 'characters' defined in this way
have the form
$\frac{\chi(g)}{\dim\chi}$, where 
$\chi(g)$ are the usual irreducible characters of the group   $G$.

Consider a collection of non-negative reals
 $\alpha_j$, $\beta_j$, $\gamma$
 satisfying the conditions
$$
\alpha_1\ge\alpha_2\ge\dots,\qquad \beta_1\ge\beta_2\ge\dots,
 \qquad \sum \alpha_j+\sum\beta_j+ \gamma=1.
$$
For such collection of parameters we define a function
 $\chi_{\alpha,\beta,\gamma}(g)$ on $S_\infty$ by the formula
$$
\chi_{\alpha,\beta,\gamma}(g)=\prod_{k=2}^\infty \Bigl(\sum_j \alpha_j^k+(-1)^{j-1}\sum_j \beta_j^k \Bigr)^{r_k(g)}
,$$
where $r_k(g)$ is the number of cycles  of length $k$ of a permutation
$g$ (actually, the product in the right hand side is finite).

\begin{theorem}
\label{th:thoma}
 {\rm(}Thoma {\rm\cite{Tho1})}
 The functions $\chi_{\alpha,\beta,\gamma}(g)$ are characters of  
 the group $S_\infty$ and all the characters of $S_\infty$ have this form.
\end{theorem}

There are two equivalent but different ways to define a representation
corresponding to a character.

The first way (which was kept in mind by Thoma): these expressions 
have the form $\tr\rho(g)$, $g\in S_\infty$, for representations 
generating Murray--von Neumann factors of type 
${II}_1$. We do not explain here what this means
(see, e.g., 
\cite{Dix}, \S5, \S13.1, \cite{GHJ}), notice that explicit
constructions of such factor-representations are contained 
in \cite{VK1} and in \cite{Ver},
notice also that the Thoma characters are limits of characters
of finite symmetric groups, see \cite{VK2}.

The second way  (which was a standpoint of the present work) was proposed by G.I.Olshanski 
in \cite{Olsh-symm}. The Thoma characters correspond to unitary representations 
of the double of the symmetric group. We pass to a discussion of the double.

\sm

{\bf\punct The bisymmetric group.%
\label{ss:bisymmetric}}  Consider a Thoma character  $\chi_{\alpha,\beta,\gamma}(g)$
and a collection of formal vectors 
$v_g$ enumerated by elements $g$ of the group $S_\infty$. 
Assign inner products for these vectors by the formula
$$
\la v_g,v_h\ra:=\chi_{\alpha,\beta,\gamma}(gh^{-1})
.$$
Since the function $\chi_{\alpha,\beta,\gamma}$ is positive-definite,
a system of vectors with such inner products can be realized in a certain
Hilbert space $\cH_{\alpha,\beta,\gamma}$. It is natural to assume that the vectors
 $v_g$ generate $\cH_{\alpha,\beta,\gamma}$. Such collection of vectors 
 is unique in the following sense: if
$v_g'\in\cH'$ is another collection of vectors with such inner products,
then there exists a unitary operator 
 $U:\cH\to \cH'$ such that $Uv_g=v_g'$ (this is the standard construction
 of a Hilbert space from a reproducing kernel, see, e.g.,  \cite{Ner-gauss}, \S7.1).

 Next, notice that for each
 $p$, $q\in S_\infty$, we have
\begin{multline*}
\la v_{pg_1q^{-1}}, v_{pg_2q^{-1}}\ra= \chi(pg_1q^{-1}(pg_2q^{-1})^{-1}\ra =\chi(pg_1g_2^{-1}p^{-1})=
\chi(g_1g_2^{-1})=
\\=
\la v_{g_1}, v_{g_2}\ra,
\end{multline*}
i.e., the transformation
 $v_g\mapsto v_{pgq^{-1}}$ does not change the matrix of inner products.
Therefore there exists a unitary operator  $T(p,q)$ in $\cH_{\alpha,\beta,\gamma}$
such that 
$$
T(p,q)\,v_g = v_{pgq^{-1}} 
$$
for all  $g\in S_\infty$. It is easy to see that 
{\it $T(p,q)$ is a unitary representation of the double of the group 
$S_\infty$,
i.e., of the group  $S_\infty\times S_\infty$}.
Next, the vector  $v_e$ is fixed by the diagonal  $\diag(S_\infty)$ of the group $S_\infty\times S_\infty$.

The explicit construction of representations of the double
 $S_\infty\times S_\infty$
 corresponding to the Thoma characters is given below in Subsection \ref{ss:chipy-rep}.

\sm

{\sc Remark.} Let $\chi$ be a character of an irreducible representation
 $\rho$ of finite or compact group  $G$ in a space $W$.
 Applying the same construction to $\chi$, we get a representation 
 of the group
 $G\times G$ in the space  $\Hom(W)$ of linear operators in the space 
$W$ defined by the formula
$$T(g_1,g_2) A=\rho(g_1^{-1}) A\rho(g_2).$$
The inner product in
 $\Hom(W)$ is given by  $\la A,B\ra=\tr AB^*$. 
 The vector fixed by the diagonal subgroup is the unit operator in
 $W$,
the system $v_g$ coincides with $\rho(g)\in \Hom(V)$.
\hfill $\boxtimes$

\sm

Recall some definitions of  representation theory.

\begin{definition} Let  $G$ be a group, $K$ be a subgroup.

{\rm 1)} An irreducible unitary representation  $\rho$ of a group  $G$
in a Hilbert space  $H$ is called 
{\rm spherical with respect to $K$} or {\rm $K$-spherical} if in $H$
there exists a unique up to a scalar factor nonzero $K$-fixed vector
$v$.

\sm

{\rm 2)} Under the same notation, 
$$
\phi(g)=\la \rho(g)v, v\ra
$$
is called  {\rm a spherical function}.

\sm

{\rm 3)} If for any irreducible unitary representation  $\rho$ of the group $G$
there exists at most one 
{\rm(}up to a scalar factor{\rm)} $K$-fixed vector, then the pair $(G,K)$ 
is called {\rm  spherical%
\footnote{For unitary representations of Lie groups
the main statement about sphericity is the following
I.M.Gelfand theorem, 1950: {\it let $G/K$ be a Riemannian symmetric space,
then the subgroup 
 $K$ is spherical in $G$}. Later there were discovered other spherical pairs,
 which are obtained from Gelfand pairs by a little enlarging of $G$ or by
a little diminishing of $K$ (on this topic there was a work \cite{Kra} 
with diverse continuations). 
\newline
There is a theory of infinite-dimensional spherical pairs that can be obtained by
limits of Gelfand (or almost Gelfand) pairs, see G.I.Olshanski's works
  \cite{Olsh-symm}, \cite{Olsh-GB}.
 Later N.I.Nessonov 
  \cite{Ness-GL} observed that spherical subgroups in infinite-dimensional groups 
  can be unexpectedly 'small',  moreover, it is possible to work with such
  spherical pairs.
    In the present paper we consider numerous spherical pairs, which have no finite-dimensional
    counterparts 
   (on classical groups, see  \cite{Ner-spheric}, on the group of diffeomorphisms of the circle, see
  \cite{Ner-diff}).
  \newline
We also note that the term 'spherical' admit other meanings.}%
  $^,$%
  \footnote{Strictly speaking, we must add once more requirement to this definition:
for any   $g\notin K$ there is a spherical representation 
  $\rho$ with $K$-spherical vector $v$, such that $\rho(g)v\ne v$.
  Let us explain what is wrong.
\newline
Consider an action of the group  $\Z$ on a countable set with trivial
stabilizers. Thus we get an embedding 
 $\Z\to\ov S_\infty$. The Lieberman theorem implies that  a
 $\Z$-fixed vector in a unitary representation of 
$\ov S_\infty$ automatically is 
$\ov S_\infty$-fixed. Therefore for all spherical pairs of the type
 $G\supset \ov S_\infty$ (this paper contains many examples of such pairs),
$\Z$-fixed vectors are fixed by the subgroup  $\ov S_\infty$.}%
.}
\end{definition}

\begin{theorem}
\label{th:dva-spheric}
 {\rm a)}  The pair $\Bigl(S_\infty\times S_\infty,\diag(S_\infty)\Bigr)$ is spherical. 
 
 \sm

{\em b)} The spherical functions on $S_\infty\times S_\infty$ are the characters 
$\chi_{\alpha,\beta,\gamma}(g_1g_2^{-1})$. 
\end{theorem}

{\sc Remark}.
It is known that for many classes of groups
(finite and compact groups, semisimple real an $p$-adic
groups, nilpotent Lie groups)
irreducible representations of the direct product
$G_1\times G_2$ are precisely tensor products  $\rho_1\otimes\rho_2$
of irreducible representations of  $G_1$ and $G_2$.
Such theorem take place if at least one group has type I
 (see \cite{Dix}, Proposition 13.1.8).
 For  $G_1=G_2=S_\infty$ it is not the case.
The Thoma representations of the group   $S_\infty\times S_\infty$
are not tensor products (except the cases 
 $\alpha_1$ or $\beta_1=1$ corresponding to one-dimensional representations).
 There is no a pathology in this phenomenon.
\hfill $\boxtimes$

\sm

{\sc Remark.} Let us restrict a Thoma representation of the double to the first factor.
Consider the weakly closed operator algebra generated by operators of the representation.
It turns out that this algebra is a 
 ${II}_1$-factor. The trace on this factor is 
$$\tr(A):=\la\rho(A) v_e,v_e\ra.$$
Recall once more that the language of factors does not used in this paper.
\hfill$\boxtimes$

Next remark. Let  $r$ be an element of the complete symmetric group
 $ \ov S_\infty$. Obviously, for all  
$g_1$, $g_2\in S_\infty$ we have
$$
\la e_{rg_1r^{-1}}, e_{rg_2r^{-1}} \ra= \la e_{g_1}, e_{g_2}\ra.
$$
This means that an operator 
 $T(r,r)$ is defined for all  $r\in \ov S_\infty$. 

Define the {\it bisymmetric group} $\bfG$ as the subgroup  in $\ov S_\infty\times \ov S_\infty$, 
consisting of pairs $(g,h)$ such that $gh^{-1}\in S_\infty$. We get the following statement.

\begin{proposition}
Each Thoma representation extends to a unitary representation of the bisymmetric
group
 $\bfG$.
\end{proposition}

It turns that the bisymmetric group is a 'good object'.
It is a type I group
(G.I.Olshanski, \cite{Olsh-symm}), its unitary representations are
classified (A.Yu.Okounkov, \cite{Oko}, see below Theorem \ref{th:OO}), 
there is a substantial harmonic analysis 
\cite{KOV} related to this group.
Also, there is the following phenomenon, which has no finite-dimensional analogs.

\sm

{\bf \punct Multiplication of double cosets
and the  train category.%
\label{ss:bisymmetric-multplicativity}} 
Let $G$ be a group,  $K_1$, $K_2$ its subgroups. A {\it double coset}
is a subset in  $G$ of the form $K_1gK_2$. By $K_1\setminus G/K_2$ we denote 
the set of double cosets, i.e., the quotient set of 
 $G$ with respect to the equivalence relation 
$g\sim k_1 g k_2$,  where $k_1\in K_1$, $k_2\in K_2$.

In the diagonal $\ov S_\infty$ of the group $\bfG$ we consider subgroups
$\bfK[\alpha]:=\ov S_\infty[\alpha]$.
Now we intend to define a product
$$
\bfK[\alpha]\setminus \bfG/\bfK[\beta]\quad \times\quad \bfK[\beta]\setminus \bfG/\bfK[\gamma]\quad \longrightarrow \quad
\bfK[\alpha]\setminus \bfG/\bfK[\gamma].
$$
For this  purpose, consider the following sequence
 $\theta_j[\beta]$ in
$\bfK[\beta]$:
\begin{equation}
\theta_j[\beta]=\begin{pmatrix}
          1_\beta&0&0&0\\
          0&0&1_j&0\\
          0&1_j&0&0\\
          0&0&0&1_\infty
         \end{pmatrix}
         .
         \label{eq:theta-j}
\end{equation}
Consider double cosets
$$
\frp\in \bfK[\alpha]\setminus \bfG/\bfK[\beta],\qquad \frq\in \bfK[\beta]\setminus \bfG/\bfK[\gamma]
$$
and their representatives
 $p\in \frp$, $q\in\frq$. Consider a sequence of double cosets
$$
\frr_j:=\bfK[\alpha] \cdot p\,\theta_j[\beta]\, q \cdot \bfK[\gamma] \quad\in\quad
 \bfK[\alpha]\setminus \bfG/\bfK[\gamma].
$$

\begin{lemma}
\label{l:stabilization}
 The sequence $\frr_j$ is eventually constant. 
\end{lemma}

Lemma can be easily proved by multiplication of 
 0-1-matrices.
 
 Now we can define the product
 $\frp\circ \frq$ as
 $$
\frp\circ \frq= \frr:=\frr_N \quad\text{for sufficiently large $N$.}
 $$
 It is easy to verify 
(also, by multiplication of matrices) that the product is  associative. 
In other words, we get a category 
 $\cS=\cS(\bfG,\bfK)$ -- the  {\it train} of the pair $(\bfG,\bfK)$.
 Its objects are non-negative integers, the sets of morphisms are
 $$
 \Mor_\cS(\alpha,\beta)=\bfK[\alpha]\setminus \bfG/\bfK[\beta].
 $$
 The category $\cS$ admits a transparent
 combinatorial description  
 (which in fact was obtained in  \cite{Olsh-symm}, see below Subsection \ref{ss:chip}).
 
 \sm
 
 {\bf\punct Multiplivativity. Representations of the train category.%
 \label{ss:train}}
 Consider a unitary representation 
 $\rho$ of the bisymmetric group $\bfG$ in a Hilbert space
$H$ (the representation is assumed to be continuous on
the subgroup  $\bfK$). 
 Denote by  $H[\alpha]$ the subspace of $\bfK[\alpha]$-fixed vectors, by $P[\alpha]$
 the operator of orthogonal projection to
 $H[\alpha]$. Consider a double coset  $\frp\in \bfK[\alpha]\setminus G/\bfK[\beta]$.
 Define the operator
 $$
 \wt\rho (\frp):H[\beta]\to H[\alpha]
 .$$
 For this purpose, we choose a representative
 $p\in\frp$ and set
 $$
 \wt\rho(\frp):=P[\alpha]\,\rho(p)\Bigr|_{H[\beta]}
 .
 $$
 
 \begin{lemma}
 \label{l:independent}
 The operator  $\wt\rho(\frp)$ does not depend on a choice
 of a representative $p\in\frp$.
 \end{lemma}
 
 {\sc Proof.} Let $v\in H[\alpha]$, $w\in H[\beta]$, let
 $h_1\in \bfK[\alpha]$, $h_2\in \bfK[\beta]$.  Then
 $$
 \la v, \rho(h_1 p h_2)w\ra= \la \rho(h_1^{-1}) v,\, \rho(p) \rho(h_2)w\ra=\la v,\rho(p) w\ra
. 
 \qquad\quad\square
 $$
 
 It turns out that we get a representation of the category
 $\cS$. 
 
 \begin{theorem}
\label{th:dva-multiplicativity} 
  {\rm(}Multiplicativity theorem{\rm)}
 For each $\alpha$, $\beta$, $\gamma$ and each 
  $\frp\in \bfK[\alpha]\setminus \bfG/\bfK[\beta]$, $\frq\in \bfK[\beta]\setminus \bfG/\bfK[\gamma]$
the following identity holds
  $$
  \wt\rho(\frp\circ \frq)=\wt \rho(\frp)\circ \,\wt \rho(\frq).
  $$
 \end{theorem}
 
 Notice also that the operators 
 $\rho(\frp)$ are  {\it contractive}
 \begin{equation}
  \|\wt\rho(\frp)\|\le 1
  .
 \end{equation}

Let us define the {\it involution}  $\frp\to \frp^*$
on the category 
$$\bfK[\alpha]\setminus \bfG/\bfK[\beta]\quad\to\quad \bfK[\beta]\setminus \bfG/\bfK[\alpha],$$
it is induced by the map  $p\mapsto p^{-1}$. Evidently,
$$
(\frp\circ\frq)^*=\frq^*\circ\frp^*.
$$
Representations of the category
 $\cS$ constructed above are  $*$-representations, i.e.,
\begin{equation}
 \wt\rho(\frp^*)= \wt\rho(\frp)^*.
\end{equation}

\begin{theorem}
\label{th:dva-approximation}
The map  $\rho\mapsto\wt\rho$ is a bijection between the set of unitary
representations of the bisymmetric group  $\bfG$  and the set of
contractive 
 $*$-representations of the category  $\cS$.
\end{theorem}

{\bf\punct The purpose of the paper.%
\label{ss:tseli}} All topics and statements discussed above
were known up to the middle of  90-s  
(\cite{Tho1}, \cite{VK1}, \cite{VK2}, \cite{Olsh-symm}, \cite{Ner-cat}).
The theory of the Thoma factor-representations and representations of the bisymmetric group
were considered as a limit of the classical representation theory. 
In \cite{Ner-cat} there was observed that a multiplication
of double cosets and a multiplivativity theorem are highly general phenomena
for infinite-dimensional groups%
\footnote{A first example of such multiplication was described
by R.S.Ismagilov in 
\cite{Ism}. Now constructions parallel to the present work
exist for infinite-dimensional classical groups, see  \cite{Olsh-GB}, \cite{Ner-spheric},
\cite{Ner-coll}, \cite{Ner-invar},
for infinite-dimensional  $p$-adic groups \cite{Ner-p}, for groups of automorphisms 
of measure spaces
 \cite{Ner-bist}, \cite{Ner-spread}, \cite{Ner-poly}.}.
The main obstacle was a description of double cosets,
which for symmetric groups was obtained in  \cite{Ner-symm}, \cite{Ner-preprint},
\cite{GN}. The main topic of the present paper is a geometric description
of various train categories in  terms of 'topological field theories'.

\sm

{\bf\punct The structure of the paper.}
\S 2 contains a proof of the Lieberman theorem.

 The main 'stencil' construction is described in   \S3.
For the group $G=S_\infty\times S_\infty\times S_\infty$ and its subgroup $K=\diag (S_\infty)$
we construct the category of double cosets. We show that morphisms
are described by two-dimensional  surfaces with checkerwise colored 
triangulations and colored edges. A product of morphisms is a gluing of surfaces.
We obtain a correspondence  (to the both sides)
between unitary representations of the group
 $G$ and representations of the category of triangulated surfaces. 
We construct a family of $K$-spherical representations of
$G$ and formulate a conjecture on a description of all
spherical representations.
In this section, proofs also are stencils. They can be applied
for a wide class of pairs $(G,K)$ without serious variations.
Therefore, in 
   \S\S3-7 we concentrate on combinatorial-geometric descriptions 
  of trains (categories of double cosets).
  
  In \S4 we return to representations of the bisymmetric group. 
  In this case the combinatorial structure is reduced to one-dimensional
  objects, which were obtained in 
   \cite{Olsh-symm}. Also, we consider other examples of pairs   $(G,K)$
whose train can be described in  terms of 'one-dimensional manifolds'.

  For the bisymmetric group Olshanski and Okounkov obtained the complete classification
  of unitary representations. Its  perception requires better  acquirements
  in representation theory than the principal text of the paper, the classification is exposed
  in Addendum to \S4.

  Some pairs
 $(G,K)$ producing two-dimensional theories are discussed in  \S 5
  (here exists a wide zoo described in  \cite{Ner-preprint}).
  
  Next, consider 
  $K$ being wreath products 
  $S_\infty\ltimes (S_m^\infty)$,  $G$ can be (for instance)
  $m$ times enlarged symmetric group. We get a category whose morphisms are described
 in terms of $m$-valent bipartite graphs. 
  
  Subsequently to bisymmetric and trisymmetric groups, we discuss 
  $n$ times product $G=S_\infty\times\dots\times S_\infty$ with the subgroup  
  $K=\diag (S_\infty)$. The corresponding train admits (at least)
  two geometric descriptions, which seems to be different.  
  One variant is given in terms of surfaces tiled by $n$-gons, we briefly
  mention it in
Subsection \ref{ss:neskolko}. The same category can be described in  terms
of $(n-1)$-dimensional pseudomanifolds, this is a content of
 \S 7.
 
In \S 8 we consider the pair consisting of infinite symmetric group $G$
and a Young subgroup  $K$. In this case, an application of the general technology
leads to two-dimensional surfaces tiled by monogons (of course there is a description of the train
in less exotic terms). This is the simplest representative of our zoo,
and  we can achieve a better understanding. We present a proof
of the formula for $K$-spherical functions on $G$ obtained by N.I.Nessonov.

\sm

{\bf\punct Remarks to \S1.} a) Different proofs of the Thoma theorem were obtained  in
\cite{Tho1}, \cite{VK2}, \cite{Oko}.

\sm

b) The problem on characters of the Thoma type can be formulated
for an arbitrary group  $G$. The interesting theory exists
for infinite unitary, orthogonal and symplectic (quaternionic unitary)
groups, see
 \cite{Voi}, \cite{VK3}. Characters can be regarded as characters of factor-representations
 or as spherical functions on the double.

\sm

c) For a discrete group $S_\infty$ it is easy to construct innumerable
families of unitary representations, see, e.g.,  \cite{Oba}
or Addition F.4 to the Russian edition of 
 \cite{Ner-cat}. It seems that attempts to extend such activity further than a mass
 production of unitary representations fail.

d)  In works \cite{Tho1}--\cite{Olsh-symm} mentioned in the beginning of paper,
representation theory of infinite symmetric groups arises in some way or another
as a limit of representation theory
of finite symmetric groups.
Almost all constructions of the present paper (except the topic with the bisymmetric group
and the Thoma representations) are not such limits.
Therefore, we get an informal problem about 'descent down': what our constructions
(for instance, various ways of 'geometrical  encoding' of symmetric groups)
can give on the level of classical theory?

\section{The Lieberman theorem}

\COUNTERS

In this section we prove the Lieberman's classification theorem
 (Theorem \ref{th:lieberman}).
For further exposition we need not in this theorem itself but 
in Theorem \ref{th:semigroup} from
\cite{Olsh-kiado} on extensions of unitary representations of  
$\ov S_\infty$ to the semigroup of all  0-1-matrices (these two theorems
are easily reduced one to another). Moreover,
below we refer to proofs of several lemmas as to 'stencil' reasonings.

\sm

{\bf\punct Equivalent definitions of continuity.%
\label{ss:admisibility}} Consider a unitary representation
$\rho$ of the group $S_\infty$
in a Hilbert space $H$. As above, let  $H[\alpha]$ denotes the subspace of all 
$S_\infty[\alpha]$-fixed vectors.

\begin{proposition} 
\label{pr:admissible}

{\rm a)} The following conditions are equivalent:

$(*)$ $\rho$ admits an extension to a continuous representation%
\footnote{Recall that a unitary representation of a topological group
is a continuous homomorphism to the group of all
unitary operators equipped with the weak operator topology
(recall that on the group of all unitary operators the weak topology coincides with the strong operator
topology).}
 of the group $\ov S_\infty$;

\sm

$(**)$ $\cup_\alpha H[\alpha]$ is dense in $H$.

\sm

{\rm b)} If $\rho$ is irreducible, then these conditions are equivalent
to the condition:

\sm

$(***)$ $\cup_\alpha H[\alpha]\ne 0$.
 
\end{proposition}

{\sc Proof.} a) The non-obvious statement is $(*)\Rightarrow (**)$.
 Let $v\in H$, $\|v\|=1$.
Denote by $W[\alpha]$ the closed convex hull of the orbit  $S_\infty[\alpha]\cdot v$.
Since the subgroups $S_\infty[\alpha]$ form a fundamental system 
of neighborhoods of unit,
 for each
$\epsilon>0$ there exists
$\alpha$ such that for each $g\in S_\infty[\alpha]$, we have  $\la gv,v\ra>1-\epsilon$;
therefore
$\|\rho(g)v-v\|<\sqrt{2\epsilon}$. Hence, for each   $w\in W[\alpha]$, we have $\|w-v\|\le \sqrt{2\epsilon}$.
On the other hand, any closed convex subset in a Hilbert space
contains a unique element of the minimal length. 
Denote such element of
 $W[\alpha]$ by $w^\circ$. By the uniqueness, 
$w^\circ$ is  $S_\infty[\alpha]$-invariant.
Thus, we can approximate 
$v$ by elements of   $\cup_\alpha H[\alpha]$ with arbitrary precision .

\sm

$(**)\Rightarrow (*)$. Let $v\in H$, $\|v\|=1$. Let $w\in  H[\beta]$ satisfy $\|w-v\|\le\delta$.
Then for all
$g\in S_\infty[\beta]$ we have $\|\rho(g)v-w\|=\|\rho(g)v-\rho(g)w\|\le\delta$, 
in particular,  $\|\rho(g)v-v\|<2\delta$.

\sm

b) Define the subgroup  $S_\alpha\subset S_\infty$ 
consisting of permutations fixing 
$\alpha+1$, $\alpha+2$, \dots.
If $\beta>\alpha$, then the subgroups
 $S_\alpha$ and $S_\infty[\beta]$ commute,
 therefore the subspace  $H[\beta]$ is $S_\alpha$-invariant. Hence 
the subspace $\cup H[\beta]$ is $S_\alpha$-invariant for each $\alpha$. Therefore
$\cup H[\beta]$ is invariant with respect to the union
of the subgroups $S_\alpha$, i.e., the whole group   $S_\infty$. By the irreducibility,
the closure of
$\cup_\alpha H[\beta]$ coincides with $H$.
\hfill $\square$

\sm

{\bf \punct Extension to the semigroup.%
\label{ss:extension-semigroup}} Denote by  $B_n$ the semigroup of
0-1-matrices of order $n$,
by $B_\infty$ the semigroup of infinite  0-1-matrices, 
equipped with the weak operator topology%
\footnote{In our case a sequence $g^{(k)}$ converges to $g$,
if for each $i$, $j$ the sequence of matrix elements  $g^{(k)}_{ij}$
coincides with $g_{ij}$ starting some place.
}. 
Notice, that the multiplication in 
 $B_\infty$ is separately continuous%
\footnote{In the space of operators with norm  $\le 1$
equipped with the weak  topology the multiplication is separately
continuous.},
but not jointly continuous.

The following statement can be easily verified
(a formal proof is in 
\cite{Ner-cat}, Lemma VIII.1.1)

\begin{lemma}
\label{l:SB-density}
 The group $S_\infty$ is dense in    $B_\infty$.
\end{lemma}

 For a matrix  $p\in B_\infty$, denote by  
 $$\{p\}_n\in B_n$$
 its left upper corner of size
 $n\times n$. The following statement is obvious.

\begin{lemma}
\label{l:SB}
The map  $p\mapsto \{p\}_n$ induces a bijection 
$$
S_\infty[k]\setminus S_\infty/S_\infty[k]\quad \to\quad B_k.
$$
\end{lemma}

Let us introduce a notation for the following element of $B_\infty$
\begin{equation}
\theta[\beta]:=\begin{pmatrix}1_\beta&0\\0&0_\infty \end{pmatrix}.
\label{eq:thetabeta}
\end{equation}
Emphasize that
\begin{equation}
\theta[\beta]=
\lim_{j\to\infty} \theta_j[\beta]
\label{eq:thetabeta-1},
\end{equation}
where $\theta_j[\beta]$ is given by the formula  (\ref{eq:theta-j}).

\begin{theorem}
\label{th:semigroup}
{\rm a)}
Each unitary representation 
 $\rho$ of the group $\ov S_\infty$ admits a unique
 extension  to a continuous representation 
 $\ov \rho$ of the semigroup
 $B_\infty$, moreover $\ov\rho(p^*)=\ov\rho(p)^*$.
 
 \sm
 
 {\rm b)} $\ov\rho\bigl(\theta[\alpha]\bigr)=P[\alpha]$, where  $P[\alpha]$ is the operator
 of orthogonal projection
 to $H[\alpha]$. 
\end{theorem}

{\sc Proof.} a)
Let $p\in B_\infty$.
Choose an element $g_n\in S_\infty$ such that  $\{g_n\}_n=\{p\}_n$.
Denote by $\frg_n\in S_\infty[n]\setminus S_\infty/S_\infty[n]$
the double coset containing  $g_n$.
By Lemma \ref{l:SB}, we have $g_{n+k}\in \frg_n$. Therefore,
\begin{equation}
P[n] \,\rho(g_{n+j})\,P[n]= P[n]\, \rho(g_{n})\,P[n].
\label{eq:PAP}
\end{equation}
(see the proof of Lemma \ref{l:independent}). Define the operators 
$$A_n(g):=P[n] \,\rho(g_{n})\,P[n].$$
Equality (\ref{eq:PAP})  takes on the form
$$P_n A_{n+j} P_n=A_n.$$
Keeping in the mind that  $\|A_n\|\le 1$, we get that the sequence
 $A_n$ has a weak limit, we denote it by   
 $\ov\rho(p)$. 

 It is easy to watch that the map
 $p\mapsto \ov\rho(p)$ is a weakly continuous map from 
 $B_\infty$ to the semigroup of operators with
 norm $\le 1$. By the separate continuity of the multiplication in
$B_\infty$ and in the semigroup of contractive operators, the continuous extension
of the representation is a representation.

\sm

b) First, $\theta[n]^2=\theta[n]=\theta[n]^*$. Therefore $\ov\rho\bigl(\theta[n]\bigr)$
is an operator of orthogonal projection.
Denote its image by $W$.

Second, in virtue of 
(\ref{eq:thetabeta-1})
 and the continuity of the representation $\ov\rho$,
the sequence
$\rho\bigl(\theta_j[n]\bigr)$ converges to $\ov\rho\bigl(\theta[\beta]\bigr)$. But
$\rho\bigl(\theta_j[n]\bigr)$ is an identical operator on  $H[n]$, hence
$H[n]\subset W$.

Third, let $g\in S_\infty[n]$, $w\in W$. Then
$$
\rho(g)w=\rho(g)\,\ov\rho\bigl(\theta[n]\bigr)\, w= 
\ov\rho\bigl(g\theta[n]\bigr) w
=\ov\rho\bigl(\theta[n]\bigr) w=w,
$$
i.e, $W\subset H[n]$. \hfill $\square$

\sm

{\bf\punct The action of $B_n$ in $H(n)$.%
\label{ss:action-B}}
For  $r\in B_n$ define the element
$\pi(r)=  \begin{pmatrix}
                      r&0\\
                      0&0_\infty
                     \end{pmatrix}\in B_\infty$.
  Denote by $B^\circ_n$ the subsemigroup   in  $B_\infty$ consisting of matrices $\pi(r)$.                 
Obviously, $\pi(r)\,\theta[n]=\theta[n]\,\pi(r)=\pi(r)$. 
Substituting these matrices to  $\ov\rho$, we get that the operator 
$\ov\rho\bigl(\pi(r)\bigr)$ has the block structure 
$$
\ov\rho\begin{pmatrix}
                      \xi_n(r)&0\\
                      0&0
                     \end{pmatrix}:H[n]\oplus H[n]^\bot\to H[n]\oplus H[n]^\bot
,$$
where $\xi_n(r)$ is an operator in $H[n]$.

\begin{corollary}
\label{cor:}
{\rm a)}
For any unitary representation 
 $\rho$ of the group   $\ov S_\infty$
 and each $n\ge 0$ we have the natural representation $\xi_n(r)$ of the semigroup $B_n$
 in the subspace  $H[n]$. 
 
 \sm
 
 {\rm b)} In $H[n]^\bot$ the semigroup $B^\circ[n]$ acts by zero operators.
\end{corollary}

\begin{lemma}
\label{l:irr-irr}
 If $\rho$ is irreducible, then $\xi_n(r)$ also is irreducible. 
\end{lemma}

{\sc Proof.} Assume that there exists a $B_n$-invariant subspace
$V\subset H[n]$. We claim that for any  
 $g\in B_\infty$ we have  $P[n]\, \ov\rho(g)\, V\subset V$. Indeed, let
 $v\in V$. Then
$$
P[n]\, \ov\rho(g)\, v=P[n]\, \ov\rho(g)\, P[n]\, v=\ov\rho\bigl(\theta[n]\bigr)\,
\ov\rho(g)\, \ov\rho\bigl(\theta[n]\bigr)\,v
=\ov\rho\bigl(\theta[n]g\theta[n]\bigr)\,v
$$
Notice that 
$\theta[n]g\theta[n]\in B_n^\circ$, and therefore $P[n]\, \ov\rho(g)\, v\in V$.
\hfill $\square$

\begin{lemma}
\label{l:determined}
 Let  $H[n]$ be
 $S_\infty$-cyclic, i.e., the system of subspaces  $\rho(g) H[n]$,
 where $g$ ranges in $S_\infty$,
 is total in $H$. Then the representation
 $\rho$ is uniquely determined by the representation $\xi_n$ of the semigroup $B_n$ in $H[n]$.
\end{lemma}

{\sc Proof.} Let $v$, $w\in H[n]$. Then
\begin{multline*}
\la \rho(g)\, v,\, \rho(h)\, w\ra= 
\la\rho(g)\,\rho(\theta[n])\, v,\, \rho(h)\,\rho(\theta[n])\, w\ra
=\\=
\la(\rho(\theta[n])\,\rho(h^*)\,\rho(g)\,\rho(\theta[n])\, v,\, w\ra=
\la\rho(\theta[n] h^* g\theta[n])\, v,\, w\ra
=\\= \la\xi_n(\theta[n] h^* g\theta[n])\, v, w\ra,
\end{multline*}
and we see that inner products 
 $\la \rho(g) v, \rho(h) w\ra$ are uniquely determined by the representation 
$\xi_n$. Hence the system of vectors  $\rho(g) v$, where $g$ ranges in $S_\infty$
and  $v$ ranges $H[n]$, is uniquely determined up to the unitary equivalence.
Therefore the operators
$$\rho(h): \,\, \rho(g) v\mapsto \rho(hg) v$$
are uniquely determined.
\hfill $\square$

\sm

{\bf\punct  The Lieberman theorem.%
\label{ss:lieberman}} Fix $n=0$, $1$, \dots\,\, and an irreducible representation $\mu$ 
of a finite group $S_n$ in the space $L=L(\mu)$. Denote by $\iota$ the trivial one-dimensional
representation of the group
$\ov S_\infty[n]$
Consider the representation of the group
 $\ov S_\infty$ induced from the representation 
$\mu\otimes \iota$ of the subgroup $S_n\times S_\infty[n]$.
Let us explain what means the term 'induced representation'%
\footnote{%
Formally, the unitary induction in the Mackey sense
(see \cite{Kir}, \S13.2)
is defined for locally compact groups.
In our case, it can be applied  
(see \cite{Kir}, \S13.2, Problems 5-7) since
$S_\infty\bigr/(S_\alpha\times S_\infty[\alpha])$
is a countable space with a discrete measure.
Also, it is possible to apply a construction with fiber
bundles
(see, e.g, \cite{Kir}, \S 13.4)
or a construction for finite groups from
 \cite{Ser}, \S7.} in this specific case.  

Consider a set 
 $\Omega$ of all $n$-element subsets in  $\N$ 
and the space  $V$ of all  $L$-valued functions on  $\Omega$.
For each  $g\in S_\infty$, $I\in \Omega$ consider a substitution   $\sigma(g,I)$ 
determined by the following rule. Denote elements of $I$ by $i_1<\dots<i_n$, 
and elements of  $gI$ by $j_1<\dots<j_n$. For each $i_m$, the element 
$gi_m$ has the form  $j_l$. We set $\sigma(g,I)m=l$.
The action of the group  $S_\infty$ in the space  $V$ we define by the formula 
$$
T_{n,\mu}(g)f(I)=\mu\bigl(\sigma(g,I))\,f(gI).
$$
We equip $V$ with the $\ell_2$-inner product,
$$
\la f_1,f_2\ra_V= \sum_\Omega \la f_1(I),f_2(I)\ra_L.
$$

\begin{theorem}
\label{th:lieberman}
{\rm a)} Each irreducible unitary representation 
of the group  $\ov S_\infty$ has the form
 $T_{n,\mu}$.
 
 \sm
 
 {\rm b)} Any unitary representation 
  $\rho$ of the group $\ov S_\infty$ is a direct sum of irreducible representations.

{\rm c)}  This decomposition is unique 
 {\rm(}i.e., collections of summands with multiplicities
 are always the same{\rm)}.
\end{theorem}

{\sc Proof.} Let us say that a {\it depth of a representation}
is the minimal $n$, for which $H[n]\ne 0$. Let the depth of 
$\rho$ be $n$. 

\sm

a) In $H[n]$ we have the action of the semigroup
 $B[n]$. Moreover, for $g\in B_n\setminus S_n$, we have $\xi_n(g)=0$ 
(in the opposite case
$n$ is not minimal). Therefore, we have in   $H[n]$ an irreducible representation
  $\mu$ of the group 
$S_n$, moreover $\rho$ is uniquely determined by a representation
 $\mu$.  Therefore,
  $\rho=T_{n,\mu}$.

\sm

b) Let us explain how to choose an irreducible subrepresentation.

First, let a representation
 $\xi_n$ of the semigroup   $B_n$ in $H[n]$ be irreducible. Consider the $B_n$-cyclic
 span  $R$ of the subspace $H[n]$.
Assume that $R$ is reducible, $R=R_1\oplus R_2$. Then the projection 
operators from $H[n]$ to $R_1$ and  $R_2$ are  $B^\circ_n$-intertwining,
and their images are subrepresentations in 
 $H[n]$, which are irreducible. Therefore,
 one of projections (say, the projection to  $R_2$) is zero. Therefore, $H[n]\subset R_1$,
 and the cyclic span of 
$H[n]$ is contained in $R_1$. We come to a contradiction.

Let $\xi_n$ be reducible, $V$ be irreducible $B_n$-subrepresentation in  $H[n]$.
Then (see the proof of Lemma  \ref{l:irr-irr})
$B_\infty$-cyclic span of the subspace 
$V$ is contained $V\oplus H^\bot$.  As we have seen above,
the cyclic span is an irreducible subrepresentation.

Repeating the same reasoning, we decompose
$\rho$ into a sum of irreducible subrepresentations. 

\sm 

c) Consider a decomposition of  $\rho$ into irreducible subrepresentations
and decompose  $\rho$ into the sum
 $\rho=\oplus_{j=0}^\infty\rho_j$, where $\rho_j$ is the sum of all irreducible summands
 of the depth  $j$. Notice that 
 $\rho_0$ is a subspace of $S_\infty$-invariant vectors.
 It is uniquely determined by the representation $\rho$, 
 and this determines  $\oplus_{j=1}^\infty\rho_j$. Next, $\rho_1$ is a cyclic span
 of the subspace 
 $H[1]$, therefore  $\rho_1$  also is determined by the representation $\rho$, etc.
 
Thus, all the summands $\rho_j$ are canonically determined. 
The problem of a decomposition of  $\rho_j$ into irreducible subrepresentations
is equivalent to a decomposition of the representation of the group 
 $S_j$ in  $H[j]$, this decomposition is unique. 
\hfill$\square$

\sm
 
{\bf \punct Remarks to \S2.}  The complete unitary, complete orthogonal, and complete symmetric groups
possess similar property \cite{Olsh-fin}, another similar group 
is the group of all measure-preserving transformations of a Lebesgue space with measure, see
 \cite{Ner-cat}, Chapter VIII.

\section{The trisymmetric group and two-dimensional triangulated bordisms}

\COUNTERS

{\bf \punct The product of double cosets.%
\label{ss:product-tri}} Consider the product 
$G=S_\infty\times S_\infty\times S_\infty$ of three copies  
of the infinite symmetric group.
We will denote elements of $G$ by 
 $\bfg=(g_{red}, g_{yellow}, g_{blue})$, the notation 
$g_\nu$ will denote one of these three 'colour' permutations.
Denote by  $K\simeq  S_\infty$ the diagonal subgroup, i.e., the subgroup
consisting of triples  $(g,g,g)$. Let $K[\alpha]$ be subgroups 
corresponding 
$S_\infty[\alpha]$ under the isomorphism
 $S_\infty\to K$. Define the elements 
 $\theta_j[\beta]\in K[\beta]$  by the same formula 
(\ref{eq:theta-j}) as above.

For each $\alpha$, $\beta$, $\gamma\in\Z_+$ we define the multiplication
of double cosets 
$$
K[\alpha]\setminus G/K[\beta]\quad\times \quad K[\beta]\setminus G/K[\gamma]
\quad \to\quad K[\alpha]\setminus G/K[\gamma]
$$
as above. Precisely, take double cosets
 $\frp$, $\frq$, choose their representatives  
$p\in \frp$, $q\in\frq$, and consider the following 
sequence of double cosets 
\begin{equation}
\frr_j=K[\alpha]\cdot p\, \theta_j[n]\, q\cdot K[\gamma]
.
\label{eq:frr-j}
\end{equation}

\begin{lemma}
\label{l:tri-correct}
 The sequence  $\frr_j$ is eventually constant.
 Its limit does not depend on the choice of the representatives $p$ and $q$.
\end{lemma}

\begin{figure}
 $$
\epsfbox{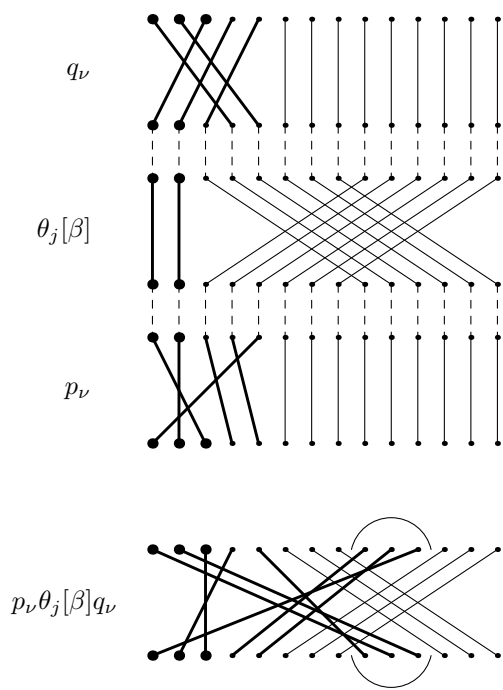}
$$
\caption{Forcing apart. To the definition of the product of double cosets.
\label{fig.1}}
\end{figure}

{\sc Proof.} Denote $I[\beta]=\{1,2\dots,\beta\}$.

We say that the {\it support} $\supp(g)$ 
of the permutation  $g$ is the set of all  
$i\in \N$ such that $gi\ne i$. 
Notice that
$\supp(g)=\supp(g^{-1})$. Denote by   $s(g)$
the maximal element of the support.
 As  sufficiently large $j$ in (\ref{eq:frr-j})
we take any number larger than  
$\max (s(p_\nu), s(q_\nu))-\beta$. Then for all  $\nu$
we have
$$
\Bigl(
\supp\bigl(\theta_j[\beta]\, p_\nu\bigr) \setminus I[\beta]
\Bigr)
\cap
\Bigl( \supp(q_\nu) \setminus I[\beta]\Bigr)=\varnothing
.
$$

Explain the phenomenon of multiplication on Fig.
\ref{fig.1}. We draw permutations 
 $q_\nu$, $\theta_j[\beta]$, $p_\nu$ and their product. 
 The vertical dashed lines denote the gluing of pictures
 corresponding  $q_\nu$, $\theta_j[\beta]$, $p_\nu$.
The fat points are elements of the sets
 $I[\alpha]$, $I[\beta]$, $I[\gamma]$ 
(on the figure, $\alpha=3$, $\beta=2$, $\gamma=3$, $j=6$).
All remaining elements of the corresponding copies of  $\N$
we call 'milky'. 'Vertical' arcs 
in  $q_\nu$, $\theta_j[\beta]$, $p_\nu$ are highlighted by fat if it least
one of its ends is fat. Also, for 
 $q_\nu$, $p_\nu$ we highlight by fat arcs whose ends are contained in the supports.
Informally, we say that fat arcs correspond to the elements of figure
containing the input information.

We call 'milky' all other vertical arcs.
Under our choice of 
 $j$ milky end of a fat arc of the permutation
 $q_\nu$ can not be joined  (through  $\theta_j[\beta]$)
with a fat arc of $p_\nu$. And vice versa.

On the lowest piece of the figure, where the permutation
$p_\nu\, \theta_j[n]\, q_\nu$ is displayed we highlight by fat
arcs that were obtained by gluing of arcs, at least one of 
which was fat. 

Change $j$ to $j+1$ and look what happens with
the fat elements of $p_\nu\, \theta_j[n]\, q_\nu$.
Milky ends of fat arcs (both upper and lower) move two unit right.
It will be remembered that  we have deal with double cosets.
We can perform a permutation of milky points of the upper row and a permutation
of milky points of the lower row and return all the ends  of the fat
arcs to the initial positions. Moreover, these permutations are the same for
all 
$\nu=1,2,3$. This proves the first statement. 

Next, consider a product $p_\nu\, h\, \theta_j[\beta]\, q_\nu$, where
$h\in K[\beta]$. 
For sufficiently large
 $j$, we have  
$$
p_\nu\, h \,\theta_j[\beta]\, q_\nu=
 p_\nu \, \theta_j[\beta]\, q_\nu (\theta_j[\beta]^{-1} \,h\, \theta_j[\beta])
, $$
and this proves the second statement.
\hfill $\square$

\begin{lemma}
\label{l:tri-ass}
The product of double cosets is associative
\end{lemma}

We omit a proof. A formal proof of Lemmas
\ref{l:tri-correct}, \ref{l:tri-ass} on the level of  0-1-matrices
is done in 
\cite{GN}.

\sm

The map $g\mapsto g^{-1}$ induces an  involution $\frg\mapsto\frg^*$ 
$$
K[\alpha]\setminus G/K[\beta] \quad \to\quad K[\beta]\setminus G/K[\alpha]
.$$
Obviously, 
$$
(\frp\circ\frq)^*=\frp^*\circ\frq^*.
$$

\begin{figure}
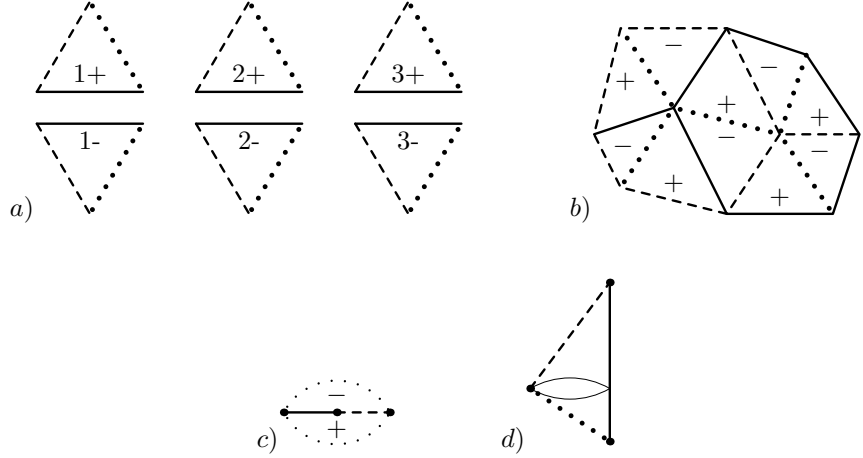

 $$
a) \epsfbox{treug.1}\qquad\qquad b) \epsfbox{treug.2}
$$

$$
c)  \epsfbox{treug.3}\qquad\qquad
d)   \epsfbox{treug.4}
$$
\caption{a) Two countable collections of colored  triangles. 
\newline
b) A result of a gluing.
\newline
c) Such position of triangles is admissible.
\newline
d) A double triangle.%
\label{fig.2}}
\end{figure}

\sm

{\bf\punct The group $G$ and combinatorial data.%
\label{ss:triangulations}} Fix an element $g=(g_{red}, g_{yellow}, g_{blue})$.
Consider a countable collection of identical triangles whose 
sides are colored red, yellow, blue clockwise, see Fig. \ref{fig.2}.
Numerate the triangles by natural numbers an mark them by the sign
 $+$.
 Consider another collection of numerated triangles whose sides are colored red, yellow, blue anticlockwise.
Mark them by the sign  $-$. The triangles are considered as oriented.

Let a permutation  $g_{red}$ send $k$ to $l$. 
Then we glue the red side of  $k$-th plus-triangle with the red side 
of   $l$-th minus-triangle according the orientations. We repeat this for 
all  $k\in \N$ and all colours (red, yellow, blue). As a result, we obtain 
a two-dimensional triangulated surface.

Note that a sense of the term '{\it triangulation}' can slightly vary.
 Specify that each triangle is embedded to the surface injectively,
 but an intersection of two triangles is not necessary a vertex or side.
Moreover, if    $m$ is sufficiently large, then  $g_\nu m=m$ for all colours,
and therefore 
$m$-th plus-triangle is glued along its perimeter 
with  $m$-th minus-triangle; thus we get a sphere separated by an 'equator; into
two triangles. We call a sphere glued from two triangles by a {\it double triangle}.

We call by an 
{\it equipped surface} $\Xi$ a union of a countable family
of oriented triangulated compact two-dimensional surfaces provided by the following combinatorial structure:

\sm

$\bullet$ each triangle is marked by  assign 'plus' or 'minus'; pluses and minuses 
are arranged checkerwise (i.e., neighbors of plus-triangles are
minus-triangles and vice versa);

\sm

$\bullet$ edges of of the triangulation are colored red, yellow, blue
in such a way that a triangle can not have edges of the same colour, moreover the three colours
of sides of plus-triangles are arranged clockwise, and colours of sides of minus-triangles
anticlockwise;

\sm

$\bullet$ all plus-triangles are numerated
by natural numbers; minus-triangles also are numerated by
the natural numbers;

\sm

$\bullet$ all but a finite number of connected components 
of the surface  $\Xi$ are double triangles such that numeric labels 
of plus-triangles and minus-triangles coincide. 

\sm

Equipped surfaces are considered up to the natural combinatorial equivalence preserving all the 
described structures.

\begin{theorem}
\label{th:tri-correspondence}
The constructed above map from the group
 $G=S_\infty\times S_\infty\times S_\infty$ to the space of all equipped surfaces
 is a bijection. 
\end{theorem}

For a proof it is sufficient to construct the inverse map. Consider an equipped
surface. For each red edge
 $\zeta$ we look at the numeric label 
 $k(\zeta)$ on the incident plus-face and to the label    $l(\zeta)$
 on the incident minus-face. Next, we set   
$g_{red}k(\zeta)=l(\zeta)$. Permutations $g_{yellow}$ and $g_{blue}$ are 
defined in the same way.
\hfill $\square$

\sm

{\sc Examples and remarks.} 1) If $\bfg$ is diagonal element, then
the surface  $\Xi$ consists of double triangles. 

2) Replacing infinity by a finite number  $n$, we get a correspondence
between the group 
$S_n\times S_n\times S_n$ and the set of equipped surfaces
glued from 
 $2n$ triangles.

3) It is easy to see that each vertex is incident to edges of only two colours and colours of these
edges interlaced, see Fig.
 \ref{fig.2}.b. Color the vertex to the complementary (to colours of edges)
 colour. It can be readily checked that red vertices
 are in one-to-one correspondence with independent cycles
  of the  permutation  $g_{yellow}^{-1} g_{blue}$.

4) Consider the subgroup  $\Gamma_\bfg$ in $S_\infty$ generated by the elements
 $g_{yellow}^{-1} g_{blue}$, $g_{blue}^{-1} g_{red}$.
 It is easy to see that components of the surface are in one-to-one correspondence 
 with orbits of the group
 $\Gamma_\bfg$.

5) Let $\bfg\in G$, $\bfh\in K$. The operation of right multiplication
 $\bfg\mapsto \bfg\bfh$ is reduced to a permutation of plus-labels on the equipped surface.
\hfill $\boxtimes$

\sm

{\bf\punct Combinatorial description of double cosets.%
\label{ss:double-cosets-surfaces}} The last remark implies the following corollary.

\begin{theorem}
\label{th:tri-cosets}
A pass from  $\bfg$ to the corresponding double coset
 $\frg\in K[\alpha]\setminus G/K[\beta]$
is equivalent to forgetting of plus-labels with numbers  $>\beta$ 
and minus-labels with numbers  $>\alpha$.
\end{theorem}

Pro forma we present a description of combinatorial data corresponding double cosets. Denote by
 $\cL[\alpha,\beta]$ the set of compact triangulated surfaces equipped with the following
 data:

$\bullet$ on each triangle we have a sign 
'plus' or 'minus', and plus-triangles and minus-triangles are arranged
checkerwise;

$\bullet$ edges of the triangulation are colored red, yellow, blue;
on plus-triangles  red, yellow, blue edges are arranged clock-wise, on minus-triangles
anticlockwise;

$\bullet$ there are fixed an injective map   $\iota^+$ from the set  
 $\{1,\dots,\beta\}$ to the set of plus-triangles and an injective map
$\iota^-$ from
$\{1,\dots,\alpha\}$ to the set of minus-triangles; we call triangles labeled in this way
 {\it entries} and {\it exits};

$\bullet$ each double triangle has at least one label. 

\sm

{\sc Examples and remarks.} 1) To keep or not to keep
the countable family of double triangles without labels 
is a test question. In any case, such triangles keep no information.

 2) Let $K$ be a group. The double coset space
 $$\diag(K)\setminus K\times K\times K/\diag(K)$$
is in one-to-one correspondence with conjugacy classes
of  $K\times K$ 
with respect to $\diag(K)$. 

3) In particular, we get a combinatorial description 
of pairs of substitutions  $\sigma_1$, $\sigma_2$ described up to a common conjugation, 
$(\sigma_1,\sigma_2)\sim (\tau \sigma_1 \tau^{-1},\tau \sigma_1 \tau^{-1})$.
\hfill $\boxtimes$

\sm

\begin{figure}
$$
\epsfbox{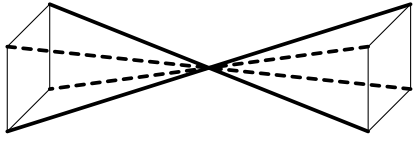}
$$
\caption{To \ref{ss:multiplication-surfaces}. A surface with glued vertices.%
\label{fig:non-normal}}
\end{figure}

{\bf \punct Multiplication of double cosets.%
\label{ss:multiplication-surfaces}} 
Let $\frP\in \cL[\alpha,\beta]$, $\frQ\in \cL[\beta,\gamma]$.
Denote a new coset  $\frR\in  \cL[\alpha,\gamma]$
according the following rule.
We make holes on the place of labeled plus-triangles 
(entries) $\iota_\frP^+(k)\in \frP$
and minus-triangles (exits) $\iota_\frQ^-(k)\in \frQ  $, 
and glue together boundaries of $\iota_\frP^+(k)$ and
$\iota_\frQ^-(k)$ according coloring of edges.  
We perform such manipulation for each  $k$ and get
a two-dimensional object equipped with a triangulation, coloring of
edges, pluses, minuses and labels on triangles.
The object consists of triangles glued by edges,
each edge is incident to a precisely two 
triangles. Such object is a two-dimensional triangulated surface,
but sum pairs of its vertices can be identified.%
\footnote{Such effect can arise if two entries having a common vertex are glued
with two exits having a common vertex}%
, See Fig. \ref{fig:non-normal}. We cut such pastings of vertices
and get a new surface with a desired structure.

\begin{proposition}
The described operation corresponds to the multiplication of double cosets.
\end{proposition}

{\sc Proof.} Consider two elements of the group  $G$ and the associated surfaces.
Glue the surfaces according our prescription.
It is easy to see that this corresponds to a multiplication
in the group.

'Evaluate' in this way the product
$$
\bfp\,\theta_j[\beta]\,\bfq = \bfp\,\bigl(\theta_j[\beta]\,\bfq\bigr)
$$
for large  $j$.
Before the gluing we must permute minus-labels on the surface corresponding to  
$\bfq$. After this,  minus-triangles of the surface  
$\theta_j[\beta]\,\bfq$ with numbers $>\beta$
 glued with double-triangles of the surface 
 $\bfp$. Under such manipulation, a surface remains to be the same,
 only labels on re-glued triangles change. 
Similarly, plus-triangles of the surface   $\frp$ with labels $>\beta$ are glued
with double triangles of the surface X $\theta_j[\beta]\,\bfq$. Notice that 
changing of label,  which took place, takes no matter,
these labels will be forgotten after a pass to double cosets.

Thus, we observe that surfaces are actually glued by triangles with labels
 $\le\beta$, and this is our operation. 
\hfill $\square$

\sm

{\sc Example.} The multiplication in the semigroup
$K[0]\setminus G/K[0]$ is an operation of disjoint union
of surfaces. In particular, this semigroup is commutative.
\hfill $\boxtimes$

\sm

The following statement is obvious.

\begin{proposition}
 The involution $\frp\to\frp^*$ corresponds to the change of the orientation of the surface
 and replacing of all minuses to pluses and all pluses to minuses. 
\end{proposition}

\sm

{\bf\punct Multiplicativity theorem.%
\label{ss:multiplicativity-tri}} 
Let $\rho$ be a unitary representation of the group  $G$ in a Hilbert space  $H$. Let $H[\alpha]$, $P[\alpha]$,
$\wt\rho(\frg)$ are introduced as above, see Subsection \ref{ss:train}. 
Then the following multiplicativity theorem holds:

\begin{theorem}
\label{th:multiplicativity-tri}
 For each $\alpha$, $\beta$, $\gamma\in \Z_+$ for each
 $$
 \frp\in K[\alpha]\setminus G/K[\beta],\quad
 \frq\in K[\beta]\setminus G/K[\gamma]
 ,$$
  we have
 \begin{equation}
 \wt\rho(\frp)  \wt\rho(\frq)= \wt\rho(\frp\circ\frq).
 \label{eq:miltiplicativity-tri}
 \end{equation}
 Moreover,
 \begin{equation}
\wt \rho(\frp^*)=\wt\rho(\frp)^*,\qquad \|\wt\rho(\frp)\|\le 1
.
 \label{eq:add-tri}
 \end{equation}
\end{theorem}

{\sc Proof.} The same arguments as in the proof of Proposition  
 \ref{pr:admissible}.b
 show that the closure 
  $V$ of the subspace  
$\cup H[\alpha]$ is $G$-invariant. Therefore the orthogonal complement  $V^\bot$ 
also is $G$-invariant.
Hence the operators   $\wt\rho(\frp)$ depend only of the restriction of
 $\rho$ to the subspace   $V$. 
Therefore, without loss of generality, we can assume that  $V=H$. 
By Proposition 
 \ref{pr:admissible},
the representation of the subgroup $K\simeq S_\infty$ in $H$ extends continuously
to the complete symmetric group 
$\ov S_\infty$. Now, applying Theorem \ref{th:semigroup}.b, we observe that the sequence
$\rho\bigl(\theta_j[\beta])$ weakly converges to the projector  $P[\beta]$. 

Choose representatives of double cosets,
$p\in\frp$, $q\in\frq$. For sufficiently large  $j$, we have 
\begin{multline*}
\ov \rho\bigl(\frp\circ\frq\bigr)\, P[\gamma]= 
P[\alpha]\,\rho\bigl(p \theta_j[\beta] q \bigr)\, P[\gamma]
=
\lim_{k\to\infty} P[\alpha] \rho\bigl(p \theta_k[\beta] q \bigr)\, P[\gamma]
=\\=
\lim_{k\to\infty} P[\alpha]\rho(p)\rho(\theta_k[\beta])\rho (q)  P[\gamma]
=  P[\alpha]\,\rho(p)\, \Bigl( \lim_{k\to\infty}\rho(\theta_k[\beta])\Bigr)\,\rho (q)
\,  P[\gamma]
=\\=
  P[\alpha]\,\rho(p)\, P[\beta]\,\rho (q)\,  P[\gamma]= 
\Bigl( P[\alpha]\,\rho(p)\, P[\beta] \Bigr)\,\Bigl( P[\beta]\,\rho (q)\, 
P[\gamma]\Bigr)
=\wt\rho(\frp)\,\wt\rho(\frq)\, P[\gamma],
\end{multline*}
the limits means the weak operator limit.
This proves (\ref{eq:miltiplicativity-tri}). The statements
(\ref{eq:add-tri}) are obvious. 
\hfill $\square$

\sm

{\bf\punct Sphericity.%
\label{ss:sphericity-tri}} 

\begin{theorem}
\label{th:tri-spheric}
 The pair $(G,K)$ is spherical.
\end{theorem}

{\sc Proof.} The semigroup $K\setminus G/ K$ is Abelian. Therefore any
its irreducible representation compatible with the involution is one-dimensional.
This semigroup acts in  $H[0]$. If its representation in  $H[0]$ is reducible then  
(see the proof of Lemma \ref{l:irr-irr}) the representation of $G$ is reducible.
\hfill $\square$

\sm

	It should be noted the following statement just mentioned and applied.

\begin{theorem}
\label{th:spheric-multiplicative}
 $K$-spherical functions of the group  $G$ are homomorphisms of the semigroup 
 $K\setminus G/K$ to the multiplicative group of complex numbers
whose absolute value 
 $\le 1$.
\end{theorem}

\sm

{\bf\punct The trisymmetric group.%
\label{ss:tri-symmetric}}
Next, note that there was a dissonance in the proof
of Theorem
 \ref{th:multiplicativity-tri}.
 The space  $H$ was decomposed into a sum of two   $G$-invariant subspaces,
$H=V\oplus V^\bot$. Concerning the representation in   $V$
we claim a nontrivial statement. Concerning 
$V^\bot$ we refer to the identity  $0\cdot 0=0$. 
Accordingly, it is natural to restrict a class
of representations, which we investigate.

We call by the 
{\it trisymmetric group} the subgroup  $\bfG$ in 
$\ov S_\infty\times\ov  S_\infty\times\ov  S_\infty$ consisting of triples
$(g_{red}$, $g_{yellow}$, $g_{blue})$ such that $g_{red} g_{yellow}^{-1}$, 
$g_{yellow}g_{blue}^{-1}\in S_\infty$. In other words, we take the subgroup 
in $\ov S_\infty\times\ov  S_\infty\times\ov  S_\infty$ generated by  
$G=S_\infty\times  S_\infty\times S_\infty$ and the diagonal subgroup  $\bfK$
consisting of elements 
$(g,g,g)$ with $g\in\ov S_\infty$. On the group  $\bfG$ we introduce the topology
assuming that the topology on 
$\bfK$ is the standard topology of $\ov S_\infty$,
and the whole group is a disconnected union of double cosets  $\bfG/\bfK$.

The following statement is tautological.

\begin{lemma}
{\rm a)} The unitary representation $\rho$ of $G$ in the space  $H$
admits a continuous extension to the group  $\bfG$ if and only if 
the restriction of $\rho$ to 
 $K$ admits a continuous extension to  $\bfK$.

 \sm
 
 {\rm b)} A unitary representation of the group $G$ in a space  $H$ admits a continuous
 extension to  $\bfG$ if and only if the subspace   $\cup H[\alpha]$ is dense $H$.
\end{lemma}

As above, we define the family of subgroups 
 $\bfK[\alpha]$ in $\bfK$. As above, the multiplication of double cosets is well-defined,
moreover, the natural map
$$
K[\alpha]\setminus G/ K[\beta]\quad \to\quad \bfK[\alpha]\setminus \bfG/ \bfK[\beta]
$$
is a bijection. 

As above, for an element of the group
$\bfG$ we draw a colored triangulated surface with labels. 
As before, it is a disjoint union of a countable family of triangulated surfaces,
and almost all surfaces are double triangles. 
The only difference is that we do not require coincidence of labels on sides of almost all
double triangles.

\sm

{\bf \punct The approximation theorem.%
\label{ss:approximation}}
Consider the category $\cK(G,K)$ -- the {\it train} of the pair $(G,K)$,
whose objects are non-negative integers, and morphisms are double cosets 
$K[\alpha]\setminus G/K[\beta]$. Consider its  $*$-representation  $R$. 
This means that for each  $\alpha\in \Z_+$ we have 
a Hilbert space 
$V[\alpha]$, and for each morphism  $\frp:\beta\to\alpha$ we have a bounded linear operator 
$R(\frp):V[\beta]\to V[\alpha]$ such that
$$
R(\frp\circ \frq)=R(\frp)R(\frq), \qquad R(\frp^*)=R(\frp)^*, \qquad R(Id_\alpha)=1_{V[\alpha]}
,$$
where $Id_\alpha$ is the unit automorphism of an object  $\alpha$, and $1_{V[\alpha]}$ 
is the unit operator in $V[\alpha]$.

\begin{theorem}
\label{th:tri-approximation}
 Any $*$-representation  $R$ of the category $\cK(G,K)$ is equivalent to
 the representation $\wt\rho$ obtain from some
unitary representation  $\rho$ of the group 
 $\bfG$, moreover, the representation  $\rho$ is unique. 
\end{theorem}

In the next subsection we obtain the representation
 $\rho$
as a limit of representations of semigroups $\Mor(\alpha,\alpha)$. 

\sm

{\bf\punct The inverse construction.%
\label{ss:inverse-construction}}
The following construction is a special case of abstract Theorem
 VIII.1.10 
from \cite{Ner-cat}.

Let $\alpha\le \beta$.
Define the following morphisms of our category:

$\bullet$ $\theta^\beta_\alpha\in \Mor(\beta,\beta)$
is the collection of double triangles with the following labels:
$$(1,1), (2,2),\dots, (\alpha,\alpha),\,\,  (\alpha+1, \varnothing),\dots, (\beta, \varnothing),\,\, 
(\varnothing,\alpha+1),\dots, (\varnothing,\beta);$$

$\bullet$ $\mu^\beta_\alpha\in \Mor(\alpha,\beta)$ is the collection of double
triangles with the following labels:
$$(1,1), (2,2),\dots, (\alpha,\alpha),  \,\, 
(\varnothing,\alpha+1),\dots, (\varnothing,\beta);$$

$\bullet$ $\nu^\beta_\alpha\in \Mor(\beta,\alpha )$ is defined as
$\nu^\beta_\alpha=\left(\mu^\beta_\alpha\right)^*.$

It is easy to see that
\begin{multline*}
\left(\theta^\beta_\alpha\right)^*= \theta^\beta_\alpha= \left(\theta^\beta_\alpha\right)^2; \\
\mu^\beta_\alpha\circ \nu^\beta_\alpha=   \theta^\beta_\alpha,\quad 
\nu^\beta_\alpha\circ  \mu^\beta_\alpha=1_\alpha;
\\
\theta^\beta_\alpha\circ 
\mu^\beta_\alpha = \mu^\beta_\alpha,
\qquad \nu^\beta_\alpha\circ  \theta^\beta_\alpha=\nu^\beta_\alpha.
\end{multline*}
Apply to these identities the functor
 $R$.  We get that 
$R(\theta^\beta_\alpha)$ is  the operator of orthogonal projection
in $V[\beta]$, and $\mu^\beta_\alpha$ is an operator
of isometric embedding 
 $V[\alpha]\to V[\beta]$,
 identifying
$V[\alpha]$ with the image of the operator $\theta^\beta_\alpha$.

Next, we construct an embedding of semigroups 
 $\zeta:\Mor(\alpha,\alpha)\to\Mor(\beta,\beta)$.
To a surface $\in\Mor(\alpha,\alpha)$ we add a collection of double
triangles with labels%
\footnote{Note that all elements of the image are not invertible.}
$$
(\alpha+1, \varnothing),\dots, (\beta, \varnothing),\,\, 
(\varnothing,\alpha+1),\dots (\varnothing,\beta).
$$
It is easy to see that
$$
\zeta(\frp)\circ  \theta^\beta_\alpha=\theta^\beta_\alpha \circ \zeta(\frp),
\qquad
\zeta(\frp)\circ \mu^\beta_\alpha=\mu^\beta_\alpha \circ \frp.
$$
Decompose the space 
$V[\beta]$ as a direct sum   $\im R(V[\alpha])\oplus \im R(V[\alpha])^\bot$.
By the above identities, the operator
 $R(\zeta\bigl(\frp)\bigr)$ has the following block structure 
$$
\begin{pmatrix}
 R(\frp)&0\\
 0&0
\end{pmatrix}
.$$
Next, notice that for
 $\alpha\le\beta\le\gamma$, we have 
$$
\mu^\gamma_\beta\circ \mu^\beta_\alpha =\mu^\gamma_\alpha.
$$
Consider the following chain of embeddings of Hilbert spaces
$$
\dots \longrightarrow V[\alpha] \longrightarrow  V[\alpha+1] \longrightarrow\dots,
,$$
take the union of all members of this chain, denote the complement by 
  $V$. In each $V[\alpha]$  we have the action of the semigroup
 $\Mor(\alpha,\alpha)$, the actions are compatible, therefore 
 in the limit space we get a representation $R^\diamond$ of the semigroup
   $\Gamma$, which is obtained as the inductive limit
 of the chain 
$$
\dots \longrightarrow \Mor(\alpha,\alpha) \longrightarrow \Mor(\alpha+1,\alpha+1) \longrightarrow\dots
.$$
By $P[\alpha]$ we denote the projection operator to   $V[\alpha]$. Then
$$
P[\alpha]\Bigr|_{V[\beta]}=R^\diamond(\theta_\alpha^\beta).
$$

\begin{lemma}
\label{l:szhim}
 For each $\frp\in\Gamma$, we have $\|R^\diamond(\frp)\|\le 1$.
\end{lemma}

{\sc Proof.} Let $\frp\in \Mor(\alpha,\alpha)$. Let the corresponding surface 
consists of $2N$ triangles. Arranging labels
on empty triangles in arbitrary way, we get 
$g\in S_N\times S_N\times S_N$ such that
$
\frp=\nu_\alpha^N g\mu_\alpha^N
.$
Therefore $R^\diamond(\frp)=R(\nu_\alpha^N) R(g) R(\mu_\alpha^N)$. 
But $R(g) R(g)^*=R(g) R(g^{-1})=1_{H[N]}$, i.e., $R(g)$ is unitary. The operator $R(\mu_\alpha^N)$
is an isometric embedding, and  $R(\nu_\alpha^N)$ is adjoint to it. 
Therefore the norm of the product 
$\le 1$.
\hfill $\square$

\begin{lemma}
Let $p\in \bfG$. Let $\frp_\alpha\in \Mor(\alpha, \alpha)$ be double coset
containing $p$, we also consider it as an element of
 $\Gamma$. Then the sequence of operators 
 $R^\diamond(\frp_\alpha)$ has a weak limit and its norm
$\le 1$.
\end{lemma} 

{\sc Proof.} This follows from the compatibility 
$$P[\alpha] R^\diamond(\frp_{\alpha+1})P[\alpha]=  R^\diamond(\frp_\alpha)$$
and the previous lemma. Decompose  $V$ as a direct sum of  subspaces 
$$V=\bigoplus_{\alpha=0}^\infty \Bigl(V[\alpha]\ominus V[\alpha-1]\Bigr)$$
and represent elements of 
 $R^\diamond(\frp_\alpha)$ as block matrices.
 All blocks whose index
 $\ge \alpha$ are zero. A pass from $R^\diamond(\frp_\alpha)$ to 
$R^\diamond(\frp_{\alpha-1})$ is zeroing out of blocks
having a subscript
 $\alpha$. The converse pass is the writing these blocks,
 moreover, the norm remains to be  $\le 1$.
\hfill $\square$

\sm

For $\frp\in G$, denote  $\rho(p):=\lim_{\alpha\to\infty}\rho(\frp_\alpha)$.
Let $p$, $q\in G$, and $r=pq$.
Let all the
 6 permutations be trivial on elements $N+1$, $N+2$, \dots.
{\it We will take only  $\alpha$, $\beta$ larger than $N$}. Then
$$
\frp_\alpha\circ\frq_\beta=(\frr)_{\min(\alpha,\beta)}.
$$
Applying to the both sides the functor
$R$ and passing to the weak limit  
$$
\lim_{\alpha\to\infty}\bigl(\lim_{\beta\to\infty}....\bigr)
,
$$
we get $\rho(p)\rho(q)=\rho(r)=\rho(pq)$.

A representation of the group 
 $G$ is constructed. Its restriction to the diagonal is continuous 
 in the topology of the complete symmetric group, and therefore it admits a continuous
 extension to the whole group
$\bfG$.
\hfill $\square$

\sm

{\bf\punct Countable tensor products of Hilbert spaces.%
\label{ss:tensors}} 
See \cite{Neu}.
Let $H_j$ be Hilbert spaces,   $\xi_j\in H_j$ be fixed unit vectors. 
In each $H_j$ choose  an orthonormal basis  $e^{(j)}_k$ 
such that $e^{(j)}_1=\xi_j$.
Consider the Hilbert space
$$
\bigotimes_{j=1}^\infty (H_j,\xi_j),
,$$
whose basis is formed by formal vectors
$$
e_{k_1}^{(1)}\otimes e_{k_2}^{(2)}\otimes  e_{k_3}^{(3)}\otimes \dots
,$$
where all the factors starting some place have the form
 $e_1^{(n)}$
(and therefore the basis is countable).

For a sequence of vectors 
$v^{(j)}\in H_j$,  we define a decomposable vector 
$$
v^{(1)}\otimes v^{(2)}\otimes v^{(3)}\otimes \dots
,$$
expanding it in the basis by opening parentheses
 (and choosing summands $e_1^{(j)}=\xi_j$ from all but a finite numbers of parentheses).
The result is contained in our space iff the products
$$
\prod_{j=1}\|v^{(j)}\|,\qquad \prod_{j=1}^\infty \la v^{(j)},\xi_j\ra
$$
converge. Under these conditions, we have  
\begin{equation}
 \la \otimes_{j=1}^\infty v^{(j)}, \otimes_{j=1}^\infty w^{(j)}\ra =\prod_{j=1}^\infty \la v^{(j)},w^{(j)}\ra_{H_j}
 \label{eq:otimes}
.\end{equation}
The construction does not depend on the choice of bases
$H_j$ (formula (\ref{eq:otimes}) can be regarded as a definition of the tensor product),
but depends on distinguished vectors  $\xi_j$. 

\sm

{\bf\punct The construction of spherical representations
of the group $\bfG$. I.%
\label{ss:spheric-tri-1}} 
Consider three Hilbert spaces, $V_{red}$, $V_{yellow}$, $V_{blue}$
(they can be finite-dimensional or infinite-dimensional).
Consider the tensor product
\begin{equation}
X:= V_{red}\otimes V_{yellow}\otimes V_{blue}
\label{eq:XVVV}
\end{equation}
and a unit vector
$\xi\in  X$. Next, consider the tensor product 
$$
H:=(X,\xi)\otimes (X,\xi)\otimes \dots= 
(V_{red}\otimes V_{yellow}\otimes V_{blue},\xi)\otimes 
(V_{red}\otimes V_{yellow}\otimes V_{blue},\xi)\otimes \dots
.$$
Also, denote
$$
\Xi:=\xi^{\otimes\infty}\in H
.$$

Define an action 
$\rho$ of the group $\bfG$ in $H$.The first copy of  $S_\infty$
interchanges  red factors $V_{red}$, the second copy interchanges yellow factors $V_{yellow}$,
the third copy the blue factors $V_{blue}$. The diagonal subgroup
 $\bfK=\ov S_\infty$ acts by permutations of factors
 $(X,\xi)$. 

\sm

{\sc Remarks.} a) We can not (except some degenerate cases)
extend the action of the red copy
of $S_\infty$ to an action of the complete symmetric group.
Indeed, a transposition of two factors  $V_{red}$ in
$$
(V_{red}\otimes V_{yellow}\otimes V_{blue},\xi)\otimes
(V_{red}\otimes V_{yellow}\otimes V_{blue},\xi)
,$$
generally speaking does not fix the vector
 $\xi\otimes\xi$, therefore the product 
$(12)(34)(56)\dots$ of elementary transpositions moves the vector $\Xi$ from $H$.

\sm

b) The vector $\Xi$ is a unique   $\bfK$-fixed vector.
Indeed, complete the vector
 $\xi$  to an orthonormal basis
 $r_1=\xi$, $r_2$, $r_3$, \dots in $X$. Let 
$h$ be a fixed vector, let $\beta_{i_1 i_2\dots}$  be coefficients of its expansion
in the basis 
$r_{i_1}\otimes r_{i_2}\otimes\dots $. Then $\beta_{\dots}$ do not change under 
permutations of underscripts. If not all   $i_k$ equal 1, 
then we get a countable number of equal coefficients.
Therefore this coefficient is 0.

c) In a similar way, $H[\alpha]$ is
$$
\underbrace{(X,\xi)\otimes\dots
\otimes (X,\xi)}_{\text{$\alpha$ times}}\otimes \xi\otimes \xi\otimes\dots
.$$

d) Denote by $U(V)$ the group of all unitary operators in a Hilbert space 
$V$. Consider in $X=U(V_{red})\times U(V_{yellow})\times U(V_{blue})$
the subgroup $Q$ fixing $\xi$. Then $Q$ acts on each factor of the product 
$X\otimes X\otimes\dots$ and therefore on the whole tensor product.
This action commutes with the representation of 
$\bfG$.
\hfill $\boxtimes$

\begin{proposition}
\label{pr:cyclic-spherical}
 The cyclic $\bfG$-span of the vector $\Xi$ is an irreducible representation of 
 $\bfG$.
\end{proposition}

{\sc Proof.} This is a variant of Lemma \ref{l:irr-irr}.

\sm

{\bf \punct Spherical functions.} Let $e_i^{red}$, $e_j^{yellow}$, $e_k^{blue}$ 
be orthonormal bases 
in $V_{red}$, $V_{yellow}$, $V_{blue}$. Let
$$
\xi=\sum \alpha_{ijk} e_i^{red}\otimes e_j^{yellow}\otimes e_k^{blue}
.
$$
Let us write an expression for the spherical function
$$
\Phi(\frg):=\la \rho(g)\Xi,\Xi\ra, \qquad g\in \frg \in \bfK\setminus \bfG/\bfK
.
$$
Consider the corresponding triangulated surface
 $\frG\in \cL[0,0]$. For each edge, we assign a basis vector
$e_i^\nu$ of the corresponding colour $\nu$ to this edge.
 We call such data  {\it assignment}. Fix an assignment $\cE$.
 For each triangle $\Delta$ we denote by  
$\iota_{red}(\Delta)$, $\iota_{yellow}(\Delta)$, $\iota_{blue}(\Delta)$ the numbers
of basis vectors  on its sides. 
 
\begin{theorem}
\label{th:formula-spheric}
\begin{multline}
 \Phi(\frg)=
 \sum_{\cE} \prod_{\text{$\Delta$ is  plus-triangle}} 
 \alpha_{\iota_{red}(\Delta) \iota_{yellow}(\Delta)\iota_{blue}(\Delta)}
 \times \\ \times
  \prod_{\text{$\Delta$ is  minus-triangle}} 
\overline{ \alpha_{\iota_{red}(\Delta) \iota_{yellow}(\Delta)\iota_{blue}(\Delta)}}
\label{eq:long}
, \end{multline}
 where $\cE$ ranges all assignments.
\end{theorem}

{\sc Proof.} Let the support of  $g$ be contained in  $\{1,\dots,N\}$.
Then it  suffices to consider the representation of  
 $S_N\times S_N\times S_N$ in $X^{\otimes N}$ and evaluate its matrix element
$\la \rho(g)\xi^{\otimes N}, \xi^{\otimes N}\ra$.

To be brief, denote the red, yellow, and blue permutations by
 $p$, $q$, $r$, the elements of bases by $x_i$, $y_j$, $z_k$.
 In this notation,
$$\Xi=
\sum \left( \alpha_{i_1 j_1 k_1} \alpha_{i_2 j_2 k_2}\dots \right) (x_{i_1}\otimes y_{j_1}\otimes z_{k_1})\otimes
(x_{i_2}\otimes y_{j_2}\otimes z_{k_2})\otimes\dots
$$
$$
\rho(g)\Xi = \sum
\left( \alpha_{i_1 j_1 k_1} \alpha_{i_2 j_2 k_2}\dots \right) (x_{i_{p(1)}}\otimes
y_{j_{q(1)}}\otimes z_{k_{r(1)}})\otimes
(x_{i_{p(2)}}\otimes y_{j_{q(2)}}\otimes z_{j_{r(2)}})\otimes\dots
$$
Consider a basis vector
  $(x_{i_1}\otimes y_{j_1}\otimes z_{k_1})\otimes
(x_{i_2}\otimes y_{j_2}\otimes z_{k_2})\dots $.
Coefficients of this vector in  $\Xi$ and $\rho(g)\Xi$ equal
respectively 
\begin{equation}
 \prod\nolimits_m \alpha_{i_m j_m k_m} , \qquad\text{and}\qquad \prod\nolimits_m 
\alpha_{i_{p^{-1}(m)} j_{q^{-1}(m)} k_{r^{-1}(m)}}
.
\label{eq:alpha-alpha}
\end{equation}
Now we write
 $e^{red}_{i_m}$ on each red edge between a plus-triangle with label 
$m$ and minus-triangle with label  $p(m)$;  in a similar way we assign basis vectors
to yellow and blue edges. We
get an assignment for the triangulated surface, 
the product 
of expressions  (\ref{eq:alpha-alpha}). This is just a summand in 
(\ref{eq:long}).
\hfill $\square$

\sm

{\bf \punct Construction of spherical representations. II.%
\label{ss:spheric-tri-2}} The previous construction
admits an extension. We can consider super-tensor products instead of the usual tensor
products.

Let $V=V^{\ov 0}\oplus V^{\ov 1}$ be a linear space
decomposed
as a sum of even and odd part.
We call such object as a 
 {\it super-space}.
 Define a tensor square of
$V$ in the usual way. But we change the operator of
transposition of summands to 
\begin{multline}
 (v^{\ov 0}\oplus v^{\ov 1})\otimes (w^{\ov 0}\oplus w^{\ov 1})= 
\bigl(v^{\ov 0}\otimes w^{\ov 0}\bigr) + \bigl(v^{\ov 0}\otimes w^{\ov 1}\bigr) +
\bigl(v^{\ov 1}\otimes w^{\ov 0}\bigr) + \bigl(v^{\ov 1}\otimes w^{\ov 1}\bigr)
\mapsto\\\mapsto
\bigl(w^{\ov 0}\otimes v^{\ov 0}\bigr) + \bigl(w^{\ov 0}\otimes v^{\ov 1}\bigr)+
\bigl(w^{\ov 1}\otimes v^{\ov 0}\bigr) - \bigl(w^{\ov 1}\otimes v^{\ov 1}\bigr).
\end{multline}
Having actions of transpositions $(12)$, $(23)$, \dots,
we can define an action of the symmetric group
 $S_n$ in a tensor power
 $V^{\otimes n}$. We can say this is in another way: transporting an add vector by an odd
 vector, we change a sign.
A tensor product equipped with such action
of the symmetric group is called a 
{\it super-tensor product}. 

Now, let $V_{red}$, $V_{yellow}$, $V_{blue}$ be super-spaces. 
We define a structure of super-space in  $X$ 
in an obvious way: 
\begin{multline*}
X^{(0)}:=\Bigl( V_{red}^{(0)}\otimes V_{yellow}^{(0)}\otimes V_{blue}^{(0)}\Bigr) \oplus 
\Bigl( V_{red}^{(1)}\otimes V_{yellow}^{(1)}\otimes V_{blue}^{(0)}\Bigr) 
\oplus\\ \oplus\Bigl(  V_{red}^{(1)}\otimes V_{yellow}^{(0)}\otimes V_{blue}^{(1)}\Bigr) 
\oplus
\Bigl( V_{red}^{(0)}\otimes V_{yellow}^{(1)}\otimes V_{blue}^{(1)}\Bigr) 
,
\end{multline*}
and $X^{(1)}$ is the orthocomplement to   $X^{(0)}$.
Choose a unit vector  $\xi\in X^{(0)}$ and consider super-tensor product 
$H=(X,\xi)\otimes (X,\xi)\otimes \dots$. The complete symmetric group   $\bfK=\ov S_\infty$
acts in 
$H$ by permutations of factors. Indeed, vectors of type  $\otimes_{j=1}^\infty v_j$,
where 
 $v_j$ are purely even or odd and $v_j=\xi$ for all but a finite number of  $j$,
form a total system in  $H$. On the other hand, 
an action of an infinite permutation on such a vector is well-defined.
Action of three copies of the symmetric group $S_\infty$ we define as above.

As above, $\otimes \xi^\infty$ is a unique   $\bfK$-fixed vector,
and its cyclic span is an irreducible representation.

\sm

{\bf\punct The Fock representation of the group
of isometries of  a Hilbert space.%
\label{ss:fock}} 
For details, see, e,g., 
 \cite{Ner-cat}.
 By
 $\Isom(V)$ we denote the group of isometries of a real Hilbert space
$V$.
Evidently, isometries have the form
 $Rv=Uv+h$, where $U$ is an orthogonal operator (we denote the group of
 all orthogonal operators by $\O(V)$), and $h\in V$.

\begin{lemma}
 Let $V$ be a real Hilbert space. Then there exists a Hilbert space 
 $F(V)$ {\rm (the Fock space)} and a total system of vectors
 $\psi_v\in F(V)$, where $v$ ranges in $V$, such that 
 $$\la \psi_v,\psi_w\ra_{F(V)}=\exp\left(- \|v-w\|^2\right/2).$$
\end{lemma}

{\sc Proof.}  Let $V=\R^n$. Consider the space  $\C^n$ with the Gaussian measure 
$\pi^{-n} e^{-\|z\|^2}$. In $L^2$ 
consider a system of vectors 
\begin{equation}
 \psi_v(z)=\exp\bigl( \la z,v\ra - \la v, v\ra/2\bigr).
 \label{eq:psi-overfull}
\end{equation}
As $F(\R^n)$, we take the closure of the linear span of vectors 
$\psi_v(z)$, $v\in\C^n$
 (in fact, $F(\R^n)$ is the subspace in  $L^2$ consisting of entire functions on  $\C^N$).

Next, consider the embeddings $\C^n\to\C^{n+1}$ and the corresponding embeddings 
$J_n$ of functional spaces 
$$
J_n f(z_1,\dots, z_n,z_{n+1})=f(z_1,\dots, z_n)
.$$
It is easy to see that $J_n$  are isometric embeddings
$$L^2\bigl(\C^n,\pi^{-n} e^{-\|z\|^2}\bigr)\to 
L^2\bigl(\C^{n+1},\pi^{-n-1} e^{-\|z\|^2}\bigr),$$
inducing embeddings $F(\R^n)\to F(\R^{n+1})$. They send vectors  
$\psi_{v_1,\dots v_n}$ to $\psi_{v_1,\dots v_n,0}$.
In limit, we get a Hilbert space with a system of vectors
 $\psi_{v_1,v_2,\dots}$, where $v_j=0$ starting some place, 
their inner products are given by the formula 
(\ref{eq:psi-overfull}).
It remains to close this system in
the limit Hilbert space.

In this way, we get the space
$F(\ell_2)$, any infinite-dimensional Hilbert space can be 
identified with
$\ell_2$.
\hfill $\square$

Let
$R$ be an isometry of the space   $V$. Consider the map of the system
$\psi_v$ to itself determined by the formula  
$$
\psi_v\to\psi_{Rv}
.
$$
Then $\la \psi(Rv), \psi(Rw)\ra = \la \psi(v), \psi(w)\ra$, 
therefore the map can be extended to a certain unitary operator  $\sigma(R)$ 
in the space $F(V)$. As a result, we get a unitary representation of the group
$\Isom(V)$.

\sm

{\sc Remark.} This construction is determined by
a slightly 'immaterial' way. In fact,  the space
$F(V)$ can be realized  as a space of holomorphic functions
on the complexification of the space  $V$; the forms of operators corresponding to
orthogonal transformations and shifts are respectively
$$
 f(z)\mapsto (Uz),\qquad f(z)\mapsto f(z+h)\exp\bigl(- \la z,h\ra_{V}\bigr)
 \qquad\qquad\boxtimes
$$

It is easy to see that this representation of the group
 $\Isom(V)$ is  
$\O(V)$-spherical, the  $\O(V)$-fixed vector is  $\psi_0$, 
the spherical function equals 
$\exp\left(-\|h\|^2/2\right)$.

\sm

{\bf \punct Constructions of spherical representations. III.%
\label{ss:spheric-tri-3}}
Our next purpose is to embed the group
 $\bfG$ to the group of isometries of a real Hilbert space
 and to obtain representations of
$\bfG$ in the Fock space by a restriction of the representation of the group $\Isom(V)$.

As $V$, we take the tensor product 
$\ell_2\otimes\ell_2\otimes \ell_2$. In this space we have the tautological representation
$\pi_1\otimes \pi_2\otimes\pi_3$ 
of the direct product $\ov S_\infty\times\ov  S_\infty \times\ov  S_\infty$
(by $\pi_j$ we denote the tautological representation 
of  $j$-th copy of $\ov  S_\infty$
in $j$-th copy of $\ell_2$).
Consider the formal expression 
\begin{equation}
\bfu:=\sum\nolimits_j e_j\otimes e_j\otimes e_j
,
\label{eq:bfu}
\end{equation}
which, of course, is not contained in the space 
 $\ell_2\otimes\ell_2\otimes \ell_2$.

 Consider the following transformation of the space
 $\ell_2\otimes\ell_2\otimes \ell_2$
\begin{multline}
R_t(g_1,g_2,g_3)v= \pi_1(g_1)\otimes \pi_2(g_2)\otimes\pi_3(g_3)v
+ t\cdot \Bigl[\pi_1(g_1)\otimes \pi_2(g_2)\otimes\pi_3(g_3)\bfu-\bfu\Bigr],
\label{eq:eee}
\end{multline}
where $t$ is a real parameter. First, we notice that the expression in the 
square brackets vanishes if  
 $g_1=g_2=g_3$. Therefore for  $(g_1,g_2,g_3)\in \bfG$ this expression is contained in 
 $\ell_2\otimes\ell_2\otimes \ell_2$. Hence affine isometric transformations 
(\ref{eq:eee}) are well defined in
$\ell_2\otimes\ell_2\otimes \ell_2$. Second, we get an action
of the group  $\bfG$
(this can be easily checked in a straightforward way,
but we can also notice that formula
 (\ref{eq:eee}) looks like an expression defining a linear transform in a space, whose
 origin is shifted to the vector
 $\bfu$).

 Restricting the Fock representation of the group 
 $\Isom(\ell_2\otimes\ell_2\otimes \ell_2)$
to the subgroup $\bfG$, we get a series of representations
of the group $\bfG$; denote these representations by  
$\upsilon^{1,1,1}_t$.

Next, we have a representation of
 $\bfG$ in $\ell_2(\bfG/\bfK)=\ell_2(G/K)$. It is spherical,
 and it is reasonable to understand it as  $\upsilon^{1,1,1}_\infty$ (in any case,
 this is so on the level of the limit of spherical functions).

 The construction can be varied. Consider the action of
$S_\infty\times S_\infty \times  S_\infty$ in $\ell_2\otimes\ell_2$, 
where the first and the second factors act in the tautological way, and the third
factor in the trivial way; i.e., in fact we get an action of  
 $S_\infty \times  S_\infty$. We consider the expression   $\sum_j e_j\otimes e_j$
 and repeat the same construction. As a result, we get a series
 of representations $ \upsilon^{1,1,0}_s$ of the group $G$.

 We can repeat this omitting the first factor or the second factor in the product
$S_\infty\times S_\infty \times  S_\infty$. Denote representations obtained in this way by
$\upsilon^{1,0,1}_p$, $\upsilon^{1,0,1}_q$.

\sm

{\bf\punct Conjectures.%
\label{ss:conjectures}} 1) Let $\Upsilon$ be a representation of 
$\bfG$ in a super-tensor product. Consider the tensor product 
\begin{equation}
\Upsilon\otimes \upsilon^{1,1,1}_t\otimes \upsilon^{1,1,0}_s\otimes 
\upsilon^{1,0,1}_p \otimes  \upsilon^{1,0,1}_q.
\end{equation}
It contains a unique 
$\bfK$-fixed vector, and its cyclic span is an irreducible representation 
of  $\bfG$. 

\sm

{\it Conjecture: all $\bfK$-spherical representations of the group  
 $\bfG$ can be obtained in this way}.

2) Consider an irreducible representation 
 $\kappa$ of the group $\ov S_\infty\times \ov S_\infty \times  \ov S_\infty$
(since $\ov S_\infty$ is a type I group, $\kappa$
is a tensor product 
$\kappa_1\otimes \kappa_2\otimes\kappa_3$ of representations  $\ov S_\infty$). 

{\it Conjecture:  any irreducible unitary representation of the group $G$ is contained
in  decomposition of some representation of the type
}
\begin{equation}
\Upsilon\otimes \upsilon^{1,1,1}_t\otimes \upsilon^{1,1,0}_s\otimes 
\upsilon^{1,0,1}_p \otimes  \upsilon^{1,0,1}_p\otimes\kappa
.\end{equation}

3) {\it Conjecture: the group  $\bfG$ has type  I, and any  unitary representation
of $\bfG$
is  can be expanded as a direct integral of irreducible
representations}.

\sm

{\bf\punct Spherical characters and self-similarity.%
\label{ss:characters}} Consider the following elements of a  semigroup $\Mor(\alpha,\alpha)$.
We take the union of   
$\alpha$ double triangles with labels  $(1,1)$, $(2,2)$, \dots, $(\alpha,\alpha)$
and an arbitrary surface without labels.
The semigroup  
$\Sigma_\alpha$ of all such morphisms is the center of the semigroup  $\Mor(\alpha,\alpha)$. 

The semigroups 
 $\Sigma_\alpha$ for different $\alpha$ are isomorphic, 
the isomorphism $\pi_\alpha^{\alpha+1}:\Sigma_{\alpha+1}\to \Sigma_{\alpha}$ 
is the forgetting of two labels 
$\alpha+1$. In notation of Subsection \ref{ss:inverse-construction},
\begin{equation*}
\pi_\alpha^{\alpha+1}(\frp)=\nu^{\alpha+1}_\alpha\circ \frp\circ\mu^{\alpha+1}_\alpha 
.
\end{equation*}

Notice that $\Sigma_0=\Mor(0,0)=K\setminus G/K$.

Let  $\rho$ be an irreducible representation of the group   $\bfG$
in a space $H$. Then the semigroup
 $\Sigma_\alpha$ acts in  
 $H[\alpha]$ by scalar operators, i.e. we get a homomorphism
 $\chi_{\alpha}$ from  $\Sigma_\alpha$ 
 to the multiplicative semigroup of complex numbers whose absolute values 
$\le 1$. For any irreducible unitary representation
 $\rho$ of  $\bfG$
$$
\wt\rho\left(\pi_\alpha^{\alpha+1}(\frp)\right)=\wt\rho\left(\nu^{\alpha+1}_\alpha\right) \wt\rho\left(\frp\right) \wt\rho\left(\mu^{\alpha+1}_\alpha\right)=\chi_{\alpha+1}(\frp)\cdot 1
$$ 
We get the following statement.

\begin{proposition}
\label{pr:sfericheskij-harakter}
 The character $\chi$ does not depend on  $\alpha$.
\end{proposition}

We call $\chi:K\setminus G/K\to \C$ by the  {\it spherical character} of an irreducible representation.

\sm

In the group 
$\bfG\subset \ov S_\infty\times \ov S_\infty\times \ov S_\infty$
we define the subgroups 
$$
\bfG[\alpha]:= \bfG\cap\left(
 \ov S_\infty[\alpha]\times \ov S_\infty[\alpha]\times \ov S_\infty[\alpha]
 \right).
$$
Clearly, the subgroups 
 $\bfG[\alpha]$ are canonically isomorphic to the group   $\bfG$.

 Consider an irreducible representation 
 $\rho$ of the group $\bfG$ in a space
$H$. Restrict it to a subgroup 
$\bfG[\alpha]$. If $\alpha$ is sufficiently large, then   $H[\alpha]\ne 0$.
In $H[\alpha]$ we have the action of the semigroup 
 $\Sigma_\alpha$ of the group  $\bfG$, it coincides with the semigroup  
$\Sigma_0$ of the group  $\bfG[\alpha]$.
It is easy to verify that 
$\bfG[\alpha]$-cyclic span of any non-zero element of 
$H[\alpha]$ is a spherical representation of   $\bfG[\alpha]$.
Moreover, {\it all spherical subrepresentations of 
the restriction  $\rho\Bigr|_{\bfG[\alpha]}$ are equivalent and their spherical functions
coincide with the spherical character of the representation   $\rho$ {\rm(}for all $\alpha${\rm)}}.

\sm

{\bf\punct Several copies of the symmetric group.%
\label{ss:neskolko}} Consider the product  of 
$n$ copies of the symmetric group  $G=S_\infty\times\dots \times S_\infty$
and the diagonal subgroup  $K$. Considerations of this section
can be repeated in this case literally. We only must consider 
 $n$-gons instead of triangles.

 Notice that we must fix a cyclic order 
 on the set of copies of the symmetric group.
A change of a cyclic order leads to an equivalent theory
but a surface corresponding to an element of the group  $G$ changes.

\sm

\begin{figure}
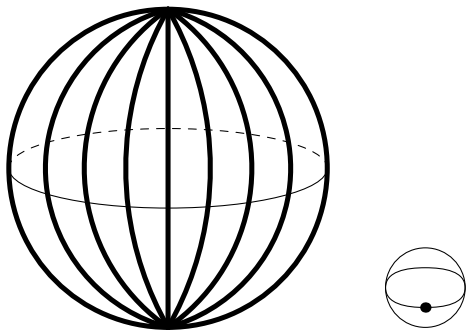

$$
\epsfbox{arbuz.1}\qquad \epsfbox{arbuz.2}
$$
\caption{Surfaces glued from digons and monogons%
\label{fig:arbuz}}
\end{figure}

We also can apply the general construction for 
$n=2$ and to get a surface glued from digons, see Fig.\ref{fig:arbuz}.
Here there is a more simple language, which is used in the next section.

Note that another language of 
 $(n-1)$-dimensional complexes is possible in a general case,
see \S7.

\sm

{\bf\punct Remarks to \S3.}
a) Main constructions of this section were obtained in
 \cite{Ner-symm}.

\sm

b) {\it The  Belyi data.} Consider a triangulated equipped surface
$\frP\in K\setminus G/K$. We claim that there is a canonical map
from $\frP$ to a double triangle  $\frT$. Precisely, we send each
plus-triangle  $\in\frP$ to the plus-triangle
of $\frT$ according the coloring of sides, we send
each minus-triangle 
 $\in\frP$ to the minus-triangle   $\in\frT$.
Further, the double triangle can be identified
with the Riemann sphere 
(vertices with the  points  0, 1, $\infty$, and edges with the segments  $[0,1]$, $[1,\infty]$, 
$[\infty,0]$). Then the complex structure can be lifted 
to the surface  $\frP$.
We get an analytic map from the Riemannian surface
 $\frP$ to the sphere $\frT$ with ramification points at vertices of the triangulation
and 3 critical values at the points   $0$, $1$, $\infty$.
According the Belyi theorem
 \cite{Bel1}, \cite{Bel2}, a map from an algebraic complex curve 
 to the sphere with 3 critical values exists precisely for curves determined
over algebraic closure  $\ov \Q$ of the field  $\Q$.

The correspondence
between $K\setminus G/K$ and the set of coverings over the
Riemann sphere with 3 critical values is a strange a posteriory coincidence of two sets.
Its meaning is unclear. Notice that by the definition the Galois group of $\ov \Q$ over $\Q$
acts on the Belyi data (i.e. on the set of coverings over the sphere with 3 critical values).
This group also remains to be invisible on our language.

 c) Consider a  group  $P$ and the space of conjugacy classes of  $P\times P$ 
by the diagonal subgroup  $P$.  The group $\GL(2,\Z)$ acts on this space 
as the group of outer automorphisms of the free group with two generators
 (see, e.g., 
 \cite{Ner-ellip}). In our case
 $P=S_\infty$, and $\GL(2,\Z)$ acts on  $K\setminus G/K$.

\section{One-dimensional constructions (chips)}

\COUNTERS

{\bf\punct $(G,K)$-pairs.}
In  \S4--\S8 we consider numerous examples of pairs 
(a group $G$, a subgroup $K$). In all the cases 
the following collection of properties take place:

\sm

1) A pair $G\supset K$ is spherical%
\footnote{Except the bisymmetric group discussed in this section,
finite analogs of all  pairs $(G,K)$ under discussion  
are not spherical. It is worth noting that a minor perturbation
of our constructions produce  large families of non-spherical pairs
$G\supset K$, see \cite{Ner-preprint}, we will not aspire to
an increasing of generality.}.

\sm

2) The category of double cosets is well defined.

\sm

3) The multiplicativity theorem holds (as Theorem \ref{th:multiplicativity-tri}).

\sm

4) An approximation theory holds (as Theorem \ref{th:tri-approximation}).

\sm

5) There is (except a degenerate case discussed in  \S8) a three-flow construction
of spherical representations as in Subsections \ref{ss:spheric-tri-1}.
\ref{ss:spheric-tri-2}, \ref{ss:spheric-tri-3}.

\sm

6) Irreducible representations admit spherical characters
 (see Subsection \ref{ss:characters}).

\sm

The argumentation of 
 \S\S2-3  can be easily applied to all the cases discussed below.
For this reason, we are concentrated on
combinatorial  realizations of double coset spaces.
The notation $G$, $K$, $\bfG$,
 $\bfK$ changes and is attributed to groups under discussion.

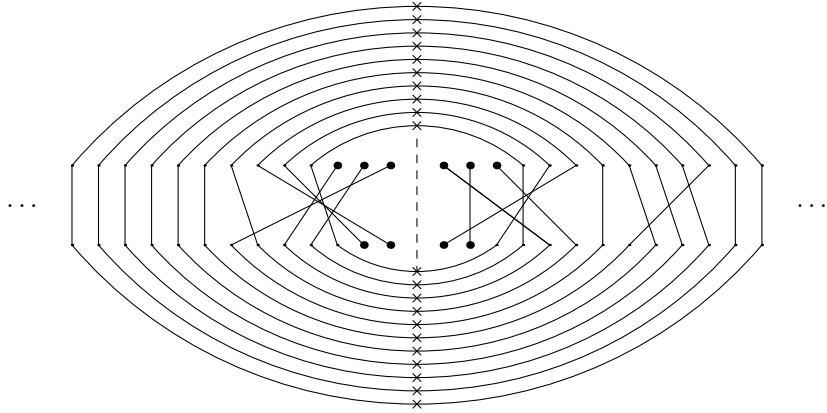
\sm

{\bf\punct Bisymmetric group and Olshanski chips.%
\label{ss:chip}}
Now, let   $G$ be the product of two copies of the symmetric group 
 $S_\infty\times S_\infty$, let $K$ be the diagonal subgroup.
 We will represent pairs of permutations as diagrams of the form
\begin{equation}
\epsfbox{bord.3}
\label{bord.3}
\end{equation}
Two copies (left and right) of natural numbers we consider as identified by
the reflection with respect to the dashed line.

\begin{figure}
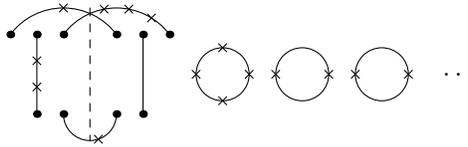

$$
\epsfbox{bord.4}
$$
a)

$$
\epsfbox{bord.5}
$$
b)

\caption{A construction of a chip from two
permutations. After gluing we get 
a collection of arcs. For each arc we remember only its 
ends (if an arc is not closed) and the number of roods on the curve.
\label{fig:chip-1}}
\end{figure}

Let us explain how to describe double cosets
$K[\alpha]\setminus G/K[\beta]$.
Consider an element  $(g_1,g_2)\in G$ and the corresponding diagram. 
For each $m>\beta$
connect by a 'horizontal' arc each left circle labeled by $m$
of the upper row with the corresponding  right circle. 
Execute the same for 
circles of the low row with numbers
 $\alpha$, see Fig. \ref{fig:chip-1},
 on the figure
 $\beta=3$, $\alpha=2$. Set a rood  on each 'horizontal' arc.
Glue horizontal arcs with 'vertical' arcs.
We call the obtained picture 
a {\it chip}. A chip can contain arcs of three types.  

a) Arcs from up to down.
Ends of such an arc are located on one side of the dashed axis,
a number of rood on the arc is even.

b) Arcs with ends  in the upper row or arc with ends in the lower row.
Ends are located in opposite sides of dashed axis, the number of roods
is odd.

c) Cycles. The number of roods is even.
Note, that the total number of cycles is infinite,
however, all but a finite number of cycles keeps precisely
two roods (and they can be removed without a loss of information).

Of course, we can replace collections of roods by non-negative integers.

\begin{figure}
$$
\epsfbox{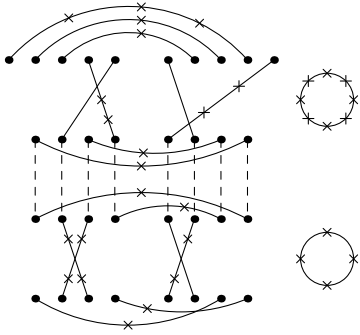}
$$
\caption{The multiplication of chips. We identify
upper circles of lower diagram with lower circles of the upper diagram.
 \label{fig:chip-multiply}}
\end{figure}

The multiplication of chips is given by gluing, see Fig.
 \ref{fig:chip-multiply}. Cycles of the factors are included to cycles of the products.

On the language of chips,
the involution
 $K[\alpha]\setminus G/K[\beta]\to K[\beta]\setminus G/K[\alpha]$ 
corresponds to the reflection with respect to the horizontal line.

\sm

{\bf \punct Representations.%
\label{ss:chipy-rep}}
Repeat the construction of Subsection
\ref{ss:spheric-tri-1}, only we must take two factors, $V$ and $W$,
instead of three in  
(\ref{eq:XVVV}). Take a unit vector  $\xi\in V\otimes W$. It is easy to verify
that there are orthonormal systems  $e_j\in V$, $e_k'\in W$ such that 
$$
\xi=\sum \alpha_j^{1/2}\, e_j\otimes e_j',\qquad \sum\alpha_j=1,\quad \alpha_j\ge 0
.
$$

If $V=V^{(0)}\oplus V^{(1)}$, $W=W^{(0)}\oplus W^{(1)}$ are super-spaces
as in Subsection \ref{ss:spheric-tri-2}, and $\xi\in V\otimes W$
is an even vector, then there are orthonormal systems 
\begin{equation} 
 e_j\in V^{(0)},\quad e_k'\in W^{(0)},\quad f_m\in V^{(1)},\quad f_l'\in W^{(1)}
 ,
 \end{equation}
 such that
$$
\xi=\sum \alpha_j^{1/2}\, e_j\otimes e_j'+
\sum \beta_m^{1/2}\, f_m\otimes f_m',\qquad \sum\alpha_j+\sum \beta_m=1,
\quad \alpha_j,\, \beta_j\ge 0.
$$
As in Subsection
 \ref{ss:spheric-tri-2} consider the super-tensor product
\begin{equation}
(V\otimes W,\xi)\otimes (V\otimes W,\xi)\otimes\dots
\label{eq:VWotimes}
\end{equation}
and the action of  $G$  in it. The vector $\Xi=\xi^{\otimes\infty}$
is a spherical vector, the numbers 
 $\alpha_j$, $\beta_m$ are the Thoma parameters (and $\gamma=0$,
 so we get not all representations).

The construction of Subsection \ref{ss:spheric-tri-3} gives
the Fock representation with
\begin{equation}
\alpha_1=e^{-t^2/2}, \,\alpha_2=\alpha_3=\dots=\beta_1=\beta_2=\dots=0.
\label{eq:toma-vyr}
\end{equation}

{\sc Remark.} A product of spherical functions in our case 
is a spherical function. Indeed, consider the tensor product
of two spherical representations of  $G$ with spherical vectors  $v$, $w$.
Then it contains a unique  
$K$-fixed vector%
\footnote{For more details,  \cite{Ner-diff}, Theorem 4.
Let
 $G\supset K$ be a spherical pair, and   $K$
 does not admit
 nontrivial finite-dimensional representations.
Then a product of spherical functions is a spherical function.
 A proof is easy.}
, and the spherical function is the product of spherical functions.
Multiplying spherical functions with parameters 
$\alpha_i$, $\beta_i$ and $\alpha_i'$, $\beta_i'$, we get the function
whose alpha-parameters are
 $\alpha_i\alpha_j'$, $\beta_i\beta_j'$, and the beta-parameters are 
 $\alpha_i\beta_j'$, $\beta_i\alpha_j'$.
\hfill $\boxtimes$

\sm

Considering tensor products of representations in super-tensor powers
 (\ref{eq:VWotimes}) with Fock representations we get all the Thoma spherical functions.
 The only exclusion is 
 $\gamma=1$. The last representation is realized in  $\ell_2(S_\infty)$,
 the group $S_\infty\times S_\infty$ acts in this space by 
 left and right shifts.

\sm

\begin{figure}
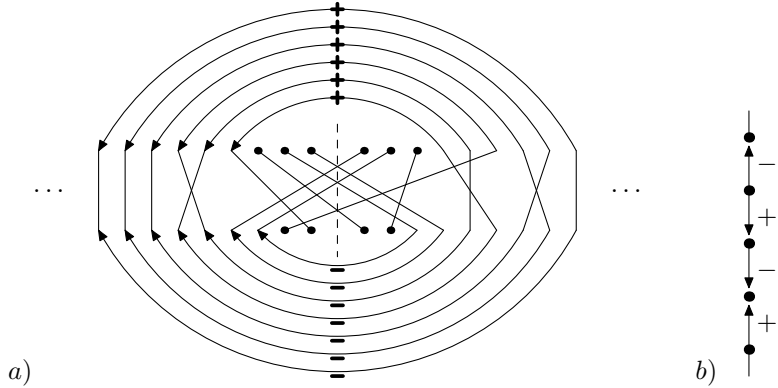

$$
a) \epsfbox{bord.7}\qquad b)\quad \epsfbox{bord.8}
$$
\caption{Reference  Subsection \ref{ss:chipy-strelki}.%
\label{fig:chipy-strelki}}
\end{figure}

{\bf\punct  Another example.%
\label{ss:chipy-strelki}} Next, consider the set $\Omega=\N\times \{0,1\}$ and
the group
 $G=S_\infty(\Omega)$. Let the subgroup 
$K$ acts on 
 $\N\times \{0,1\}$ 
 by permutations of the  factor $\N$
(as compared to the previous example,
we have the same subgroup  $K$,  and the group  $G$
is enlarged).
 We draw elements of the group  $G$ as diagrams of the form 
  \begin{equation}
  \epsfbox{bord.17} 
  \label{bord.17}
  \end{equation}
(now we allow arcs going from the left part of the picture to
the right side and vice versa).
To pass to double cosets 
(see Fig. \ref{fig:chipy-strelki}.a)
 $K[\alpha]\setminus G/K[\beta]$,
we as above connect the  elements 
of the upper row with the same  numbers  $>\beta$ by an arc.
However, now we also draw an arrow from the  right to  left
and set the sign
 $+$ to arrow. Also we draw arrow corresponding elements of 
 lower row and put to arrows signs   $-$.  

After gluing we get a collection of chains. These chains are closed
or 
are both ends are black circle on the figure.
Chains consist of interlacing links of two types, arrows with signs 
and simply segments. Segments can be contracted to points without a loss
of information (the only exception are chains consisting of a unique segment).
A piece of a chain is shown on Fig.  \ref{fig:chipy-strelki}.b. 
Notice that pluses and minuses on chains interlace. 

\sm

{\sc Remark.} In this construction, pluses and minuses are
essential only for the distinction of orientations
of cycles. 
\hfill $\boxtimes$

\sm

A rule of multiplication of double cosets is clear.

Constructions of representations can be easily extended to this case.
For instance, consider a Hilbert space  
$V$ and a unit vector $\xi\in V\otimes V$. Consider the tensor product
$$
(V\otimes V,\xi)\otimes (V\otimes V,\xi)\otimes\dots
$$
The group 
 $G$ acts by permutations of factors  $V$,
 the subgroup $K$
by permutations of blocks $(V\otimes V,\xi)$.
Action of
 $K$ admits an extension to an action
 of the complete symmetric group  $\bfK=S_\infty$.

\sm

\begin{figure}
$$
\epsfbox{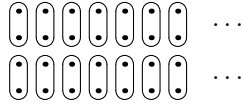}
$$
\caption{To Subsection \ref{ss:chip-3}.  The group $K_1$ acts by permutations
of columns of the upper array,
the group $K_2$ permutes columns of the lower array.\label{fig:chip-color}}
\end{figure}

\begin{figure}
$$
\epsfbox{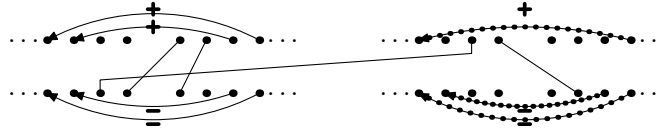}
$$
\caption{To Subsection \ref{ss:chip-3}. Constructing
of a chip. In the upper row we have 4 copies of 
$\N$ (i.e., the set  $\N\times\{1,2,3,4\}$ stretched to a line).
In the lower row we draw another copy of  the set
  $\N\times\{1,2,3,4\}$. The element  $G$ is displayed as a collection
  of strings from up to down. 
\label{fig:chip-3}
}
\end{figure}

{\bf \punct A more complicated example.%
\label{ss:chip-3}} Now consider the union  $\Omega$ of 4 copies of   $\N$,
it is convenient to think that  $\Omega=\N\times\{1,2,3,4\}$.
 Let $G=S_\infty(\Omega)$.
Consider the subgroup  $K_1\subset G$, which fixes all elements
of the set 
$\N\times\{3,4\}$, on the product   $\N\times\{1,2\}$ it acts by permutations of 
the factor $\N$. The subgroup  $K_2$ acts in the similar way  on $\N\times\{3,4\}$,
fixing $\N\times\{1,2\}$ . We display elements of this group
by  pictures of type  (\ref{bord.3}),  (\ref{bord.17}), but now 
we must draw 4 copies of $\N$. To be definite we 
arrange them in the order 
$$
\dots\,5_1\, 4_1\, 3_1 \,2_1\,1_1\quad 1_2\, 2_2\, 3_2 \,4_2\,5_2\,\dots
\qquad
\dots\,5_3\, 4_3\, 3_3 \,2_3\,1_3\quad 1_4\, 2_4\, 3_4 \,4_4\,5_4\,\dots
$$

Denote
$$
K[\alpha_1,\alpha_2]:=K_1[\alpha_1]\times K_2[\alpha_2].
$$
Thus we get a multiplication
\begin{multline*}
K[\alpha_1,\alpha_2]\setminus G/K[\beta_1,\beta_2]\,\,\,\times\,\,\,
K[\beta_1,\beta_2]\setminus G/K[\gamma_1,\gamma_2]
\,\,\,\to\\ \to \,\,\,
K[\alpha_1,\alpha_2]\setminus G/K[\gamma_1,\gamma_2]
.\end{multline*}

Let us explain a picture corresponding to a double coset, see Fig.
\ref{fig:chip-3}. We draw arrows connecting the corresponding elements
of $\N\times\{2\}$ and  $\N\times\{1\}$; color arrows, say, blue. On upper arrows 
we draw 
 $+$, on lower arrows   $-$. Execute the same operation with arrows 
 from  $\N\times\{4\}$ to $\N\times\{3\}$, we only color arrows red
(on the figure we mark 'red' by lines dashed with fat dots). Further, we glue 
chains as above.

It remains to explain how to produce representations. Consider a Hilbert space
$V$ and two unit vectors 
$
\xi\in V\otimes V$, $\eta\in V\otimes V.
$
Consider the tensor product
$$
(V\otimes V,\xi)\otimes (V\otimes V,\eta)\otimes (V\otimes V,\xi)\otimes (V\otimes V,\eta)\otimes\dots
$$
The group 
 $G$ acts by permutations of the factors  $V$, the subgroup $K_1$ 
 by permutations of factors 
 $(V\otimes V,\xi)$, and the group $K_2$ permutes factors $(V\otimes V,\eta)$.
The action of $K_1$ (resp. $K_2$) can be extended to the action
of the complete symmetric group.

 \sm
 
\begin{figure}
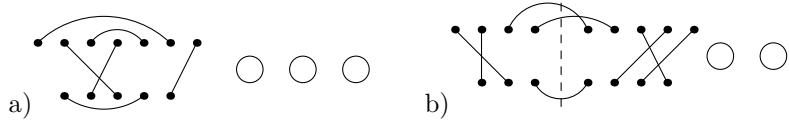

$${\mathrm a)} \epsfbox{bord.16}
\qquad
{\mathrm b)} \epsfbox{bord.18}
 $$
 \caption{Reference  Subsection \ref{ss:brauer}.a. Brauer diagrams.
\label{fig:brauer}}
 \end{figure}

 {\bf\punct Remarks to \S4.\label{ss:brauer}} a) {\it The Brauer category.} See \cite{Bra},
 see, also, \cite{Ker}.
Objects of the Brauer category are non-negative integers. 
 A morphism from $m$ to  $n$ is a diagram of the type presented on Fig.\ref{fig:brauer}.a
 with $m$ fat dots in the upper row and  $n$  fat dots in the lower row.
Products of morphisms 
  $m\to n$ and $n\to k$ corresponds to the gluing of diagrams.
 
   Consider the space
  $V=\C^p$ equipped with a symmetric nondegenerate bilinear form 
 (and respectively by the action of the orthogonal group   $\O(p)$).
By diagram we construct an operator  
 $V^{\otimes m}\to V^{\otimes n}$. For vertical arcs we assign a permutation
 of factors, for upper horizontal arcs we assign the convolutions
 with respect to the corresponding indices, to lower horizontal arcs 
 the multiplications by invariants. For a cycle we assign a multiplication by
 the scalar  $p$. As a result, we get an action of the Brauer category
 by 
 $\O(p)$-intertwining operators in tensor powers of 
 $V$.
 
 Literally the same construction takes place for symplectic groups
(in this case, to a cycle we assign the multiplication by  $(-p)$).

A similar construction exists for
 for complete linear groups. In the last case,
we take mixed tensors
 $V^{\otimes m}\otimes (V')^{\otimes k}$
 with fixed $m-k$. The corresponding set of fat points is split  into two
 parts
  ($m$ and  $k$ points). Horizontal arcs can connect 
  only points from different groups
(convolutions are possible only for pairs 
 $V$, $V'$), and vertical arcs connect points from similar
 groups, 
 see Fig. \ref{fig:brauer}. 

\sm

b) The construction of spherical representations of the bisymmetric group
arose as a result of several transformations 
of the Vershik and Kerov construction \cite{VK1}. 
The action in the space 
(\ref{eq:VWotimes}) was introduced in  \cite{Olsh-symm}.
A general construction was proposed in  \cite{DN},
we follow to a simple construction from  \cite{Ner-thoma}.

\bigskip
\begin{center}
{\bf\large Addendum to Section  4.
The classification of

representations of the bisymmetric group.}
\end{center}

\medskip

Return to notation of Subsections
\ref{ss:chip}--\ref{ss:chipy-rep}. 

\sm

{\bf\punct The duality.} Denote by  $\U(V)$ and $\U(W)$ 
the groups of all unitary operators in  $V$ and
$W$ respectively. Denote by  $\cU$ the subgroup  in $\U(V)\times \U(W)$,
consisting of all operators
 $A\otimes B$ fixing the vector
 $\xi$. It easy to see that this group
 is a direct product 
of unitary groups
$$
\cU:=\prod_{x\in [-1,1], x\ne 0} \U(n_x)
,$$ 
where $n_x$ is the multiplicity of the entry of  
$x$ into the collection $\alpha_i$, $\beta_j$. In fact,
this product is finite or countable
 (and almost all nontrivial factors have the form   $\U(1)$), 
 the group  $\cU$ is compact and separable. 
 Its irreducible representations are
 tensor products of representations of factors
 (in fact, these products are finite, see  \cite{HR}, 27.43).
 Next, we notice that pairs of operators 
$(e^{i\theta},e^{-i\theta})$ are contained in the center of 
 $\cU$ and act on  $V\otimes W$ trivially. 

The group $\cU$ acts on each factor of the product
 (\ref{eq:VWotimes}).
Since $\cU$ fixes  $\xi$,
it acts on the whole tensor product.
Obviously $\cU$ commutes with the action
of the bisymmetric group
 $\bfG$. The following counterpart of the Schur--Weyl
 duality holds, see
 \cite{Olsh-symm}.

\begin{theorem} 
{\rm {a)}}
The groups  $\bfG$ and $\cU$ are dual in the tensor product
{\rm(\ref{eq:VWotimes})}
in the following sense.
The representation of  $G\otimes \cU$ is a direct sum of 
representations 
$\rho_\nu\otimes \pi_\nu$, where $\rho_j$ are unitary representations 
of the group  $G$ and 
$\pi_\mu$ are unitary representations of $\cU$,
moreover  $\rho_\nu$ are pairwise distinct and  $\pi_\mu$ are pairwise distinct.

\sm

{\rm {b)}} The spherical representation of
$\bfG$ is realized in the space of 
$\cU$-fixed vectors.
\end{theorem}

{\bf\punct The classification of unitary representations
of the bisymmetric group.%
\label{ss:oo}} 
Let $\rho$ be a unitary representation of the bisymmetric group
 $\bfG$ in a space $H$, let
$H[n]$ be the subspaces of  $\bfK[n]$-fixed vectors. 
The representation  $\rho$
has a spherical character (see \ref{ss:characters}), 
which in our case
must be one of the Thoma characters 
(therefore the numbers $\alpha_i$, $\beta_j$, $\gamma$ are defined).

\begin{figure}
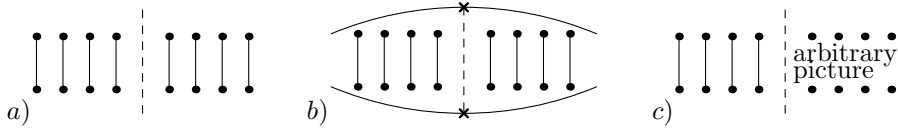

$$
a)\epsfbox{bord.11}\qquad b) \epsfbox{bord.12} \qquad c) \epsfbox{bord.13}
$$
\caption{ a) The unit chip.
b) The projector $H[n+1]\to H[n]$.
c) An element of  $\Gamma_1$.
\label{fig:for-oo}}
\end{figure}

 In subspaces $H[n]$, we have irreducible $*$-representations
 of the chip semigroups  $\bfK[n]\setminus\bfG/\bfK[n]$.
The reasonings of Lemma  \ref{l:determined} 
shows that any such representation 
uniquely determines the initial representation  $\rho$.
Consider the minimal 
$n$, for which   $H[n]\ne 0$. It is easy to verify that the chip
displayed on 
 Fig. \ref{fig:for-oo}.b, corresponds to the projection operator
$H[n]$ to $H[n-1]$. Since this operator is zero,
all the operators corresponding to chips with horizontal arcs
are zero.
Denote the semigroup of chips without
horizontal arcs by  $\Delta_n$. It  naturally splits in to a product
of three semigroups,
$$
\Delta_n= \Xi\times \Gamma_1\times \Gamma_2
,$$
here $\Xi$ is the free Abelian semigroup, whose elements are the unit chip 
 (see. Fig.\ref{fig:for-oo}.a) together with a collection of cycles.
The semigroups
 $\Gamma_1$ and $\Gamma_2$ are isomorphic; $\Gamma_1$ consists of diagrams of
 the form exposed on Fig.\ref{fig:for-oo}.c,
 where in the left-hand side there are no roods on arc; $\Gamma_2$
 consists of bilaterally reflected diagrams.
 
 The semigroup $\Xi$ is isomorphic to 
  $\bfK[0]\setminus\bfG/\bfK[0]$. Its action is determined by
  the spherical character of the representation.
 
 The semigroups $\Gamma_1$, $\Gamma_2$ are isomorphic to
 the semidirect product  $\Gamma:=S_n\ltimes \Z_+^n$.

 We denote Young diagrams by
 $\lambda$, $\mu$, \dots. The number of cells of   
  $\lambda$  we denote by  $|\lambda|$, the number of rows by $l(\lambda)$,
   the number of columns by
 $h(\lambda)$. Recall that irreducible representations
 of the group  $S_k$ have canonical enumeration by   $k$-cell 
 Young diagrams, see, e.g.,  \cite{Kir}, 16.2.

\begin{figure}
$$
\epsfbox{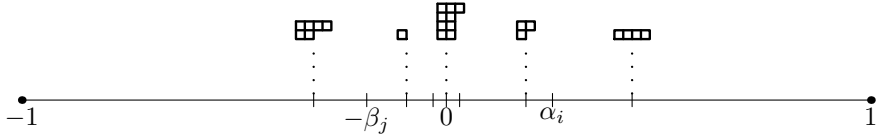}
$$
\caption{Reference  Subsection \ref{ss:oo}.
A diagram distribution.
\label{fig:oo}}
\end{figure}

We say that a {\it diagram distribution} $\Lambda$
is a finite subset   $\frl$  (the support of a distribution)  of the segment $[-1,1]$,
and to each point $x\in \frl$ we assign a Young diagram
 $\lambda(x)$.
We also require  $\sum_{x\in \frl}|\lambda(x)|=n$.

\begin{lemma}
Irreducible $*$-representations of the semigroup
 $\Gamma$ are enumerated by diagram distributions. 
\end{lemma}

This can be easily proved by the Wigner--Mackey argumentation, see
 \cite{Kir}, 13.3. Let us describe 
 the representation corresponding to a diagram
 distribution 
  $\Lambda$.
   Consider the following subsemigroup in  $\Gamma$
$$
\Gamma_\Lambda:=\prod_{x\in \frl}S_{|\lambda(x)|}\ltimes \Z_+^{|\lambda(x)|}
  .$$
  Consider the irreducible representation of a group
 $S_{|\lambda(x)|}$,
corresponding to the diagram  $\lambda(x)$, extend it to the semigroup
   $S_{|\lambda(x)|}\ltimes \Z_+^{|\lambda(x)|}$ assuming that
   the generators of the semigroup 
   $\Z_+^{|\lambda(x)|}$ act as multiplications by   $x$.
   Consider the tensor product of such  representations and induce%
\footnote{Formally, the induction is defined 
for representations of groups, and not for semigroups.
However, a definition from  \cite{Ser}, \S 7, is appropriate in our case.} 
a representation of  $\Gamma$ from the representation of the semigroup
$\Gamma_\Lambda$.

Thus an irreducible representation of the semigroup
 $\Delta_n$ is determined by the Thoma parameters $\alpha_i$,
$\beta_j$ and two diagram distributions,
say,  $\Lambda$ and $M$.
However, not all such representations are appropriate for us.

\begin{theorem}
\label{th:OO}
{\rm  (Olshanski--Okounkov, \cite{Olsh-symm},\cite{Oko})}
A representation  $\Delta_n$ determined by parameters 
$\alpha$,
$\beta$ and $\Lambda$, $M$ can be realized in 
a space  $H[n]$ of fixed vectors of an irreducible representation
of the bisymmetric group
 $\bfG$ if and only if the following conditions hold:

\sm

$\bullet$ The supports of the distributions  $\Lambda$ and $M$
are contained in the set 
$$\{\alpha_1,\alpha_2,\dots, -\beta_1,-\beta_2,\dots,0\}.$$

\sm

$\bullet$ For $x>0$ the number $l(\lambda(x))+l(\mu(x))$ does not exceed 
the multiplicity of the entry
$x$ in the sequence $\alpha_i$.
 
\sm 
 
$\bullet$ For $x<0$ the number $h(\lambda(x))+h(\mu(x))$ does not exceed
the multiplicity of the entry
$x$ in the sequence  $(-\beta_j)$.
\end{theorem}

Notice that a reconstructions of 
 $\rho$ from such data is constructive in a certain sense
 (see the  proof of Lemma  \ref{l:determined}),
 however, a question about more explicit realizations is still open.
As well, for spherical representations of  $\bfG$ with $\gamma>0$
an explicit description of the space of representation
remains to be unknown. 


\section{Two-dimensional constructions}

\COUNTERS

{\bf\punct Example 1.%
\label{ss:iznutri}}
Consider the set $\Omega=\N\times \{1,2,3\}$.
We regard elements of sets 
$\N\times \{1\}$, $\N\times \{2\}$, $\N\times \{3\}$ colored
red, yellow, blue respectively. Therefore
any element of $\Omega$ is uniquely determined by 
its number and its colour. Let  $G=S_\infty(\Omega)$.
 By $K$ we denote the group of all finite permutations
 of the set 
$\N\times \{1,2,3\}$ induced by permutations of 
the set  $\N$.

\sm

{\sc Remark.} As compared with the trisymmetric group
from
 \S3, we have enlarged the group  $G$, the group $K$ remains the same.
\hfill $\boxtimes$

\sm

Now to each element of the group 
 $G$ we will assign a triangulated surface  with a certain additional data.
 For this purpose, we consider a collection of equal triangles, whose sides
 are colored from inside,
\begin{figure}
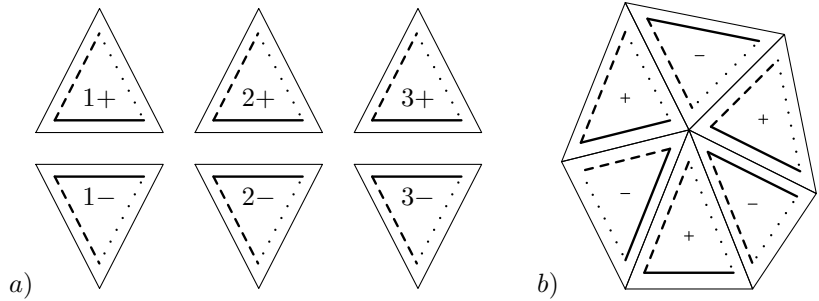

$$
a) \epsfbox{surface.1}
\qquad
b) \epsfbox{surface.2}
$$
\caption{Reference  Subsection \ref{ss:iznutri}. a) Triangles colored from inside.
 b) The surface obtained by gluing.%
\label{fig:iznutri}}
\end{figure}
in red, yellow, blue clockwise. We numerate these triangles
and set 'plus' to each of them, see Fig. \ref{fig:iznutri}.a.
Next, we draw a countable collection of reflected colored triangles, we also numerate them
and put 'minus' to each triangle. Let  $g\in G$ send an element
of  $\Omega$ with number $k$ and colour $\nu$ to an element
with number $l$ and colour  $\mu$. 
Then we glue  $\nu$-edge of $k$-th plus-triangle
and  $\mu$-edge $l$-th minus-triangle according 
the orientation. As in
 \S3, we get a triangulated  surface such that:

\sm

$\bullet$ plus-triangles and minus-triangles are arranged  checkerwise;
 
\sm 

$\bullet$ plus-triangles are numbered by $\N$; minus-triangles also;

\sm
 
$\bullet$ on interior sides of plus-triangles colors are arranged clockwise;
on interior sides of minus-triangles anticlockwise;
colours on two  sides of an edge can be different;
 
\sm 
 
$\bullet$  all components of the surface are compact;
all but a finite number of are double triangles such that numbers on 
faces coincide and colours on two sides of each  edge coincide;

\sm 
 
 It is easy to see that
  an element
 $g$ can be reconstructed from such data.
We must examine all the edges
and look to colors and labels on two sides of an edge.
 
 \sm
 
 {\it The passing to double cosets}
 $K[\alpha]\setminus G/K[\beta]$ is the forgetting of
 plus-labels  $>\beta$ and minus-labels  $>\alpha$.

 All that has been said about trisymmetric group extends to our case more-or-less
 literally.

\sm

{\bf\punct Degree of generality.%
\label{ss:generality}} Consider a finite collection of countable sets
$\Omega_1$, \dots, $\Omega_p$. Let
$$
G=S_\infty(\Omega_1)\times \dots\times S_\infty(\Omega_p)
.$$
 Let $K_1$, \dots, $K_q$ be copies of infinite symmetric group.
 Let each group $K_j$ act on each $\Omega_i$ in such a way that
 orbits are one-point sets or equivalent to  $S_\infty/S_{\infty}[1]\simeq \N$.
 Let nontrivial orbits of different groups  $K_j$ be disjoint.
 Then, on  $\Omega_1\cup \dots \cup  \Omega_p$
we get an action of the group $K_1\times\dots\times K_q$. 
In this way, we get an embedding  
$$
K_1\times\dots\times K_q \,\, \to \,\,G.
$$
Let (although this restriction  is superfluous)  
the group
$K_1\times\dots\times K_q$ has no fixed points
on  $\Omega_1\cup \dots \cup  \Omega_p$. 

Then for each element of
 $G$ we can construct a two-dimensional surface tiled
 by polygons and colored in a certain way (as in examples discussed above).
Different types of polygons also must be colored in different
colors%
\footnote{Since colors of polygons distinguishing different groups 
 $K_j$ have an origin different from colors of edges, we can introduce another attributes
 as 'perfume', 'melody', 'citizenship'.}
 corresponding to different factors $K_j$, the number of sides
 of a polygon corresponds to the number of nontrivial orbits of the group
 $K_j$.

 We omit an exposition of a formal rule, it is contained in 
 \cite{Ner-preprint}.

 All that has been said above on trisymmetric group in
 \S3 can be extended to this general situation. 

\sm

{\sc Remark.} If each group  $K_j$ has $\le 2$ nontrivial orbits, we
glue surfaces from digons, in this case the language of chips can be applied.
There exists also the case when each group  $K_j$
has a unique nontrivial orbit. This case is examined in   \S8.

\sm

\begin{figure}
$$
\epsfbox{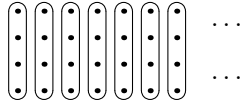}
$$
\caption{Reference  Subsection \ref{ss:wreath}.\label{fig:wreath}}
\end{figure}

{\bf \punct Wreath products.\label{ss:wreath}} Let  $H$ be a group.
Consider the group 
$H^\N$ consisting on functions  $f$  on $\N$ with values in $H$
and equal the unit for almost all 
$n\in \N$. The product is defined pointwise. The group $S_\infty$
acts on  $H^\N$ by permutations  of
$\N$. Therefore we can define 
the semi-direct product $S_\infty\ltimes H^\N$. This group is called the 
{\it wreath product} of the groups $S_\infty$ and  $H$.
We denote it by 
$$
S_\infty \wre H:=S_\infty\ltimes H^\N.
$$
In fact, below we meet the cases
$H=\Z_k$ and 
$H=S_k$. Describe the corresponding
wreath products in a more convenient form.

Consider the set
 $\N\times \{1,\dots,k\}$.
 We can regard it as the set of points of the integer lattice in a horizontal strip.
We  divide this strip into vertical columns,
see Fig.\ref{fig:wreath}.  The wreath product $S_\infty\wre S_k$
consists of permutations 
preserving the partition into columns
(i.e., we can permute columns and then take a permutation in each column).

Wreath product $S_\infty\wre \Z_k$ is the group of permutations preserving the partition
into columns and preserving cyclic order inside columns.

We introduce also the {\it completed wreath product} $\ov S_\infty\ov{\wre} H$,
which is the semi-direct product
of the complete symmetric group and the group of
all functions
 $\N\to H$.

\sm

{\bf\punct An example with wreath product.} Let  $\Omega=\N\times \{1,2,3\}$,
 $G=S_\infty(\Omega)$. As  a subgroup $K$ we
take wreath product $S_\infty\wre\Z_3$. By $K[\alpha]$ we denote the subgroup in
 $K$ fixing all elements of the first  $\alpha$ columns.

For  $g\in G$, we construct a two-dimensional surface in literally the same way as in
Subsection
\ref{ss:iznutri}. A difference appears, when we pass to double cosets.
In this case, removing label from a triangle we also must remove
the coloring of the interior side of its perimeter.

\sm

{\sc Example.} Element of $K\setminus G/K$ is a triangulated surface,
on which we arrange pluses and minuses checkerwise.
\hfill $\square$

\sm

Explain what is changing in constructions of representations.
Consider the Hilbert space 
$V$ and a vector $\xi\in V\otimes V\otimes V$, which is fixed under cyclic permutations
of factors. Next, consider the tensor product 
$$
(V\otimes V\otimes V,\xi)\otimes (V\otimes V\otimes V,\xi)\otimes\dots
$$
The group
 $G$ acts by permutations of factors  $V$, and the group $K$ by permutations of blocks
$(V\otimes V\otimes V,\xi)$ and cyclic permutations in each block
 $V\otimes V\otimes V$. The action of  $K$
 can be extend to the action of the completed wreath product
 $\ov S_\infty\ov{\wre} \Z_3$.

\section{Categories of bipartite graphs}

\COUNTERS

Here we consider only one example 
(but this construction has numerous variations).

  \sm
  
\begin{figure}
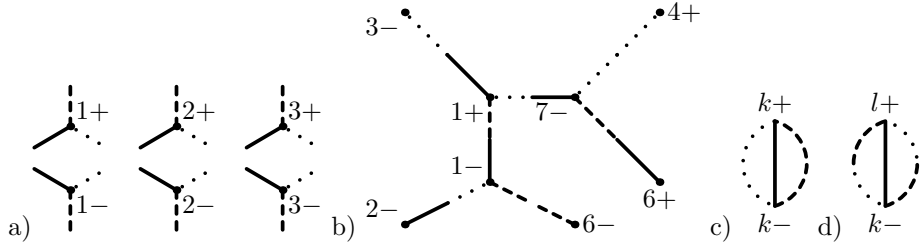

 a)\epsfbox{surface.3}
 b)
 \epsfbox{surface.4}
 c) \epsfbox{surface.5}
 d) \epsfbox{surface.6}
 \caption{Reference  Subsection  \ref{ss:bipartie}.
 \label{fig:bipartie}}
 \end{figure}

 {\bf\punct The correspondence between permutations
 and bipartite 3-valent graphs.%
 \label{ss:bipartie}}
Let $G$  be the group of finitely supported permutations of
 $\Omega:=\N\times\{1,2,3\}$, let $K$ be the wreath product of $S_\infty$ with $S_3$.
As above, we think that elements of
the sets  $\N\times\{1\}$, $\N\times\{2\}$, $\N\times\{3\}$
are colored red, yellow, blue respectively. 
We draw a countable collection of such pictures: a vertex incident to 3 segments
(we call them 
{\it semi-edges}), segments are colored red, yellow, and blue respectively.
To each vertex we set 'plus' and the number   $\in \N$.
Draw another countable family of the same pictures,
 to vertices we assign 'minuses' and 
numbers (see Fig. \ref{fig:bipartie}.a).
Let an element 
 $g\in G$ send an element of colour $\nu$ and number $k$ to an element of colour $\mu$
and number  $l$. Then we glue the semi-edge of colour  $\nu$
incident to the plus-vertex with number  $k$ to the semi-edge
of the colour  $\mu$ incident the minus-vertex with number  $l$. 
We do this for each element of  $\Omega$
and get a bipartite graph, see Fig.
 \ref{fig:bipartie}.b with the following structures and properties:

$\bullet$  on vertices we have 'pluses' and 'minuses',
moreover, edges can connect only plus-vertices with-minus-vertices; 

$\bullet$ each vertex is incident precisely 3 edges;

$\bullet$ semi-edges are colored red, yellow, blue;
semi-edges adjacent to a vertex are colored different colors;

$\bullet$ plus-vertices are numbered; minus-vertices also;

$\bullet$ the graph consists of finite components,
moreover, almost all components are
two-vertex with coinciding labels on the both vertices and coinciding colors
on halfs of each edge (see Fig. \ref{fig:bipartie}.c).
 
 It is easy to see that an element of the symmetric group
 can be reconstructed from these picture.
 
 Consider a larger group
 $\bfK\supset K$ being the completed wreath product 
 $\ov S_\infty\ov\wre S_3$
 and the group of permutations
 $\bfG$ generated by  $\bfK$ and $G$.
 To this group we can apply the same procedure of constructing of a graph.
 The list of properties of the graph obtained in this way is almost the same, 
 but the last condition must be weaken:

$\bullet$ the graph consists of finite components
and almost all components are-two-vertex
(see Fig. \ref{fig:bipartie}.d).

\sm

{\bf \punct Category of double cosets.%
\label{ss:cosets-graphs}} Let $K[\alpha]$ be the subgroup in
$K$ fixing all points of the set 
$$
\{1,\dots,\alpha\}\times\{1,2,3\}\subset \N\times\{1,2,3\}
.$$
The pass to double cosets corresponds to forgetting plus-labels with numbers 
$>\beta$ and colours of semi-edges incident to these vertices,
and also forgetting  minus-labels with numbers
$>\alpha$ and the colours of the incident semi-edges.

\sm

{\sc Example.} Passing to  $K[0]\setminus G/ K[0]$ means forgetting all the labels 
and colours. Consider the same construction for finite groups,
$G=S_{3N}$, $K=S_N\wre S_3$.
We get that 
 $$(S_N\wre S_3)\setminus S_{3N}/(S_N\wre S_3^{N})$$ 
is in one-to-one correspondence with the set of all  $2N$-vertex
bipartite three-valent graphs 
 (the plus-partite and the minus-partite remains to be distinguished). \hfill $\boxtimes$

\sm

\begin{figure}
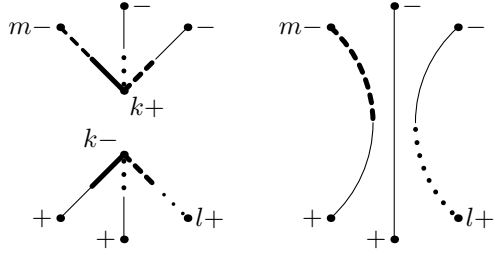

$$
 \epsfbox{surface.7}\qquad  \epsfbox{surface.8}
$$
\caption{Reference  Subsection\ref{ss:cosets-graphs}.
Gluing of three-valent graphs.
\label{fig:skleika-graphs}}
\end{figure}

Consider double cosets
$$
\frp\in K[\alpha]\setminus G/K[\beta],\quad \frq\in K[\beta]\setminus G/K[\gamma]
.$$
Describe their product on the language of graphs
 (see Fig. \ref{fig:skleika-graphs}).
Let
 $k\le \beta$, take a vertex of $\frp$ with label $k$ and minus-vertex of  $\frq$
with label  $k$. Cut from the graph these vertices together with the incident semi-edges
and glue 3 short-cut edges of the graph
 $\frp$ with 3 short-edges of short-cut of  $\frq$
 according the colours of removed semi-edges.
Executing this for all   $k$, we get a morphism  from $\gamma$  to $\alpha$.

\sm

{\bf\punct Representations of the group  $\bfG$} can be constructed as above.
For instance, we can consider (super)-Hilbert space  
 $V$, (even) vector 
$\xi\in V\otimes V\otimes V$ invariant with respect to permutations
of the factors and the action of the group 
$\bfG$ in the tensor product 
$$
(V\otimes V\otimes V,\xi)\otimes (V\otimes V\otimes V,\xi)\otimes \dots
.$$

 To repeat the construction of representations in the Fock space, we consider
a Hilbert space  
$L$ with an orthonormal basis   $e_{k,\nu}$, where $\nu=1$, $2$, $3$, 
 and $k$ ranges in  $\N$. The group $\bfG$ acts in  $L$ in the natural way.
Next, we consider the tensor product  $L\otimes L\otimes L$,
  as
 $\bfu$ (see \ref{eq:bfu}) we take
 $$
 \bfu:=\sum_{k=1}^\infty \sum_{\sigma\in S_3} 
 e_{k,\sigma(1)}\otimes  e_{k,\sigma(2)}\otimes  e_{k,\sigma(3)}
 $$
 and repeat the construction of Subsection \ref{ss:spheric-tri-3}.

 \sm
 
 {\bf\punct Variations.%
 \label{ss:var}} A) Let  $G$ be the group of all finitely supported permutations
 of the set 
 $\Omega:=\N\times\{1,\dots, l\}$, let  $K$ be the wreath product
  $S_\infty \wre S_l$. Then all the above considerations
 can be repeated in or case,
 we only get 
 bipartite  $l$-valent graphs instead of 3-valent. 
 
 B) Let
 $$G=S_\infty(\N\times\{1,\dots, 6\}),\quad 
 K=\bigl(S_\infty\wre S_3\bigr)
 \, \times\, \bigl(S_\infty\wre S_3\bigr)
 .$$
 Then we introduce two 'colours' (or 'perfumes') for coloring of vertices of a graph
 and repeat the same construction.

 C) Here it is possible a wider generality as in Subsection \ref{ss:generality}.

\begin{figure}
 $$
 \epsfbox{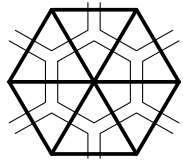}
 $$
 \caption{Reference  Subsection \ref{ss:var}. A triangulation and the dual ribbon graph.
 \label{fig:fat}}
 \end{figure}

 D) Constructions of \S5 can be exposed on the language of ribbon graphs.
 Recall that a ribbon graph is a graph with fixed cyclic order
 of edges at each vertex. There is a standard way to assign a ribbon graph
 to  any triangulated oriented surface 
(see Fig. \ref{fig:fat}).
 
 \section{Bordisms of pseudomanifolds}
 
 \COUNTERS
 
 {\bf \punct Pseudomanifolds.}
 Recall a definition of a pseudomanifold. 

 First, we determine a
 {\it simplicial cell complex}. 
 Consider a disjoint union 
 $\coprod \Xi_j$ of a finite collection of simplices 
$\Xi_j$. Consider a quotient-space   $\Sigma$ of the space 
$\coprod \Xi_j$ by a collection of equivalences
satisfying the following properties:

\sm

a) For any simplex  $\Xi_i$ the tautological map
 $\xi_i:\Xi_i\to \Sigma$ 
is an embedding. Therefore we can regard  $\Xi_i$ as a subset in $\Sigma$.

\sm

b) For any pair of simplices   $\Xi_i$, $\Xi_j$  their intersection
$\xi_i^{-1}\bigl(\xi_i(\Xi_i)\cap \xi_j(\Xi_j)\bigr)\subset \Xi_i$ 
is a union of faces  $\Xi_i$ and partially defined 
map
$$
\Xi_i\stackrel{\xi_i}{\longrightarrow} \Sigma\stackrel{\xi_j^{-1}}{\longrightarrow} \Xi_j
$$
is affine on each face.

Such objects are called 
  {\it simplicial cell complexes%
 \footnote{The standard definition of a simplicial complex
 contains a more rigid requirement:
 intersection of two simplices is empty or is a face.}.}
 
  A {\it pseudomanifold} (for details, see  \cite{ST}, \cite{Gai1})
of  dimension   $n$ is a simplicial cell complex satisfying the following
conditions:

\sm

a) Each face is contained in an  $n$-dimensional face. We call  $n$-dimensional faces
by
  {\it chambers.}

\sm

b) Each $(n-1)$-dimensional  face is contained in  precisely
two chambers
\footnote{Usually, also 'strong connectedness' is required:
the complex is connected and remains to be connected after a removing of all 
faces of codimension 
$\ge2$.}.

\sm

{\sc Remark.} Recall the origin of notion of a pseudomanifold.
 As it is well-known (R.Thom) that generally speaking a cycle
 in a manifold can not be realized as an image
 of a manifold. On the other hand, any singular $\Z$-cycle (or $\Z_2$-cycle)
 can be realized as a pseudomanifold (by the definition of a cycle, see
 \cite{ST}). Also,   intersection homologies  were introduced by Goretski
 and MacPherson through pseudomanifolds, \cite{GM}.
\hfill $\boxtimes$

\sm

{\it Link.} Let $\Sigma$  be a pseudomanifold,
$\Gamma$ its  $k$-dimensional face. Consider all 
$(k+1)$-dimensional faces  $\Phi_j$,
containing $\Gamma$ and choose a point   $\phi_j$ in
interior of each face 
  $\Phi_j$. For any face $\Psi_k\supset \Gamma$
consider the convex hull of all points  $\phi_j$,
that are contained in
 $\Psi_k$. The link of the face
$\Gamma$ is the simplicial cell complex, whose faces are such 
convex hulls. 

A pseudomanifold is called
{\it normal}, if the link of any face of codimension
 $\ge 2$ is connected.
 
A {\it normalization}
 $\wt\Sigma$ of a pseudomanifold  $\Sigma$ 
(\cite{GM})
is a normal  pseudomanifold 
 $\wt\Sigma$ and a map $\wt\Sigma \to \Sigma$ such that

\sm

--- the restriction of $\pi$ to each face   $\wt\Sigma$
is an affine bijection of faces;

\sm

--- the map  $\pi$ send  $n$-dimensional and  $(n-1)$-dimensional
faces to different faces. 
\begin{figure} 
$$
\epsfbox{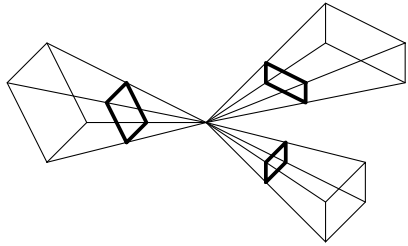}
$$ 
\caption{To the definition of the normalization. A (zero-dimensional) face
with a disconnected link.
\label{fig:link}}
\end{figure}
 
{\it Any pseudomanifold has a unique up to a natural equivalence normalization} (\cite{GM}).

Let us describe its construction. Let $\Sigma$ be non-normal. Let  $\Xi$ 
be the face of codimension 
 $\ge 2$ whose link consists of  
 $m>1$ connected components  (see Fig. \ref{fig:link}). 
Consider a small neighborhood 
 $\cO$ of the face  $\Xi$ of the pseudomanifold $\Sigma$. Then the set
 $\cO\setminus \Xi$
is disconnected and consists of  
 $m$ components, denote them by $\cO_1$, \dots, $\cO_m$.
Let  $\ov\cO_j$ be the closure of  $\cO_j$ in $\Sigma$, $\ov\cO_j=\cO_j\cup \Xi_j$.
We replace  $\cO$ by the disjoint union of sets
 $\ov\cO_j$ and get a new pseudomanifold 
$\Sigma'$. Repeat the same operation with another
face with a disconnected link.
These operations enlarge the number of faces of codimension
  $\ge 2$.
The faces of dimensions  $n$ and $(n-1)$ remains to be the same
 (and their inclusions are preserved).
Therefore, the process is finite. As a result, we get a normal pseudomanifold.

 \sm
 
 {\bf\punct Collections of permutations and pseudomanifolds.}
 As a group $G$ we take the product of $(n+1)$ copies of the group $S_\infty$,
 and  $K$ is the diagonal subgroup%
 \footnote{Recall that this pair  (group, subgroup) was considered above
 in 
 Subsection \ref{ss:neskolko}}.
 Take  $n+1$ kinds of dyes, red, orange, 
 \dots,
 blue, paint by these colors copies of the group
 $S_\infty$. Take a countable collection of
 identical  $n$-dimensional simplices, whose faces are dyed  these colours,  
 (different faces have different colours). Separately, we dye
 vertices in such a way that a colour of a vertex coincides with the colour of the opposite face.
Numerate  these simplices, and put to each simplex the sign  $+$.
Draw another collection of simplices, which differs 
from the simplices of the first collection 
by a reflection%
 \footnote{Note that all colorings of  $(n-1)$-faces of a chamber into  $(n+1)$ pairwise different
 colors are equivalent with respect to affine transformations of  $\R^n$.
 If we consider only proper affine transformations, then
 we get two non-equivalent types of simplices.}. We numerate these
 simplices and put sign <<$-$>>.
 
 Let an element of,say, green copy of
 $S_\infty$ send $k$ to $l$. Then we glue together the green face of  
 $k$-th plus-simplex with green face of     $l$-th minus-simplex
 in such a way that each vertex is identified with the vertex of the same colour.
 We execute this for all the copies of the symmetric group and all numbers
 $k$. As a result, we get a simplicial cell complex, which is 
 a countable disjoint union of 
 {\it  oriented normal}  pseudomanifolds.
 The complex satisfies the conditions:
 
 $\bullet$ chambers are marked by pluses and minuses checkerwise, i.e., a neighborhood
 (through $(n-1)$-dimensional wall) of a plus-chamber is a minus-chamber and vice versa;

 $\bullet$ plus-chambers (respectively, minus-chambers) are numerated by
 natural numbers;
 
 $\bullet$ $(n-1)$-faces are dyed to abovesaid colors. Colours of $(n-1)$-walls of 
 each chamber are pairwise different, and coloring of walls correspond to the orientation.
 Vertices also are dyed and a colour of a vertex differs from colours of all adjacent $(n-1)$-faces.

 $\bullet$ almost all components are 
 obtained by gluing of two simplices with coinciding labels.

 \sm
 
 It is easy to see that such collections of data are in one-to-one correspondence with the
 group
 $G$.
 
 \sm
 
{\bf \punct Bordisms of pseudomanifolds.} 
A pass to double cosets is  usual. We forget plus-labels with numbers 
 $>\beta$ and minus-labels with numbers   $>\alpha$.
 
 Let $\frp\in K[\alpha]\setminus G/K[\beta]$, $\frq\in K[\beta]\setminus G/K[\gamma]$.
 For any $k\le\beta$ we take the minus-chamber of the pseudomanifold
 $\frq$ and the plus-chamber of the pseudomanifold   $\frp$.
 Remove these chambers from pseudomanifolds and glue the
 boundaries of the holes according colours of  $(n-1)$-dimensional faces.  
 Performing this for all 
  $k\le\beta$, we get a new pseudomanifold with colored 
  $(n-1)$-faces and vertices, with pluses and minuses arranged 
  checkerwise and plus-labels 
 $\le \gamma$ and minus-labels  $\le\alpha$.
 Generally, this pseudomanifold is not normal.
 Its normalization corresponds to the product of double cosets
 $\frp\circ\frq$.
 
 \sm
 
{\bf\punct The correspondence with two-dimensional construction.}
Above in Subsection \ref{ss:neskolko} to the same double cosets $K[\alpha]\setminus G/K[\beta]$
we assigned two-dimensional surfaces tiled by $n$-gons.
Let us explain how to obtain the same surface starting from
the pseudomanifold. Let us numerate  $(n-1)$-faces $A_1$, \dots, $A_n$
of our model simplex $\Delta$. Choose a point inside each intersection
 $A_1\cap A_2$, $A_2\cap A_3$, \dots, $A_n\cap A_1$
 and connect these points by the closed polygonal line.
 Dye each segment into the colour of
 the $(n-1)$-face containing it. 
 Consider a two-dimensional surface $\Pi\subset \Delta$ homeomorphic
 to a disk whose boundary coincides with our polygonal line and interior 
 is contained in the interior of the simplex.
It remains to notice that all the chambers 
of our pseudomanifold are copies of the simplex $\Delta$,
and therefore in each chamber we have a copy of the surface
 $\Pi$. Uniting these surfaces, we get a two-dimensional surface tiled by 
 $n$-gons.

\sm
  
 {\bf\punct Remarks to \S7.} Constructions of pseudomanifolds from
 colored graphs were considered in
\cite{Pez},
 \cite{Fer}, \cite{FGG}. In \cite{Gai1}, \cite{BuG}, \cite{Gai2} there was used the correspondence
 between pseudomanifolds and collections of involutions. Here we follow  \cite{GN}.

\section{Spherical functions with respect to the Young subgroup. Nessonov theorem}

\COUNTERS

The object described below on the general language of
\S5 and \cite{Ner-preprint} corresponds to two-dimensional surfaces
glued from monogons, see Fig.  
\ref{fig:arbuz}.

\sm

{\bf\punct The $(G,K)$-pair with the Young subgroup.} Now, let 
$$G=S_\infty(\N\times \{1,\dots,m\}).$$
 Let $K_j=S_\infty(\N\times \{j\})$,
$\bfK_j=\ov S_\infty(\N\times \{j\})$,
and
$$K=K_1\times\dots\times K_m,\qquad \bfK=\bfK_1\times\dots\times \bfK_m.$$
Let  $\bfG$ be the group generated by   $G$ and $\bfK$.

For $g\in\bfG$, we denote by  $s_{\nu\mu}(g)$ the number of elements of a color
$\nu$ that are send by our permutation 
 $g\in \bfG$ to elements of a color $\mu$. For any  $\mu$ these numbers satisfy the relation
 \begin{equation}
 \sum\nolimits_{\nu:\nu\ne \mu} s_{\nu\mu}=\sum\nolimits_{\nu:\nu\ne \mu} s_{\mu\nu}
 \label{eq:s-mu-nu}
 \end{equation}
 (the number of elements exiting from a set $\N\times\{\mu\}$ equals the number of coming elements).

 Consider a collection of representations of the group
 $\bfG$. Consider a Euclidean space
 $V$ and $m$ unit vectors    $\xi_1$,\dots, $\xi_m\in V$ generating 
 $V$.
 Consider the tensor product
$$
 \Bigl[(V,\xi_1)\otimes (V,\xi_1)\otimes\dots]\otimes
 \Bigl[(V,\xi_2)\otimes (V,\xi_2)\otimes\dots]\otimes\dots
 .
$$
 The group $G$ acts by permutations  of factors $V$, each 
 $\bfK_j$ act by permutations in each group.
  The vector
  $$\xi_1^{\otimes \infty}\otimes\dots \otimes \xi_m^{\otimes \infty}
  $$
  is 
 $\bfK$-fixed. Its cyclic span is a  
  $\bfK$-spherical representation of  $G$. The spherical function equals  
  $$
  \Phi(g)=\prod\nolimits_{\nu,\mu:\,\nu\ne \mu} \la \xi_i,\xi_j\ra^{s_{\nu\mu}(g)}
  .
  $$
  
{\bf\punct The Nessonov theorem.}  
  
  \begin{theorem}
  \label{th:ness}
  {\rm(\cite{Ness})}
  The construction described above exhaust all
   $\bfK$-spherical representations of the group 
 $\bfG$.
  \end{theorem}
 
 A proof is given below in Subsections
 \ref{ss:char}--\ref{ss:proof}.
  
  
\sm  
  
\begin{figure}
$$
\epsfbox{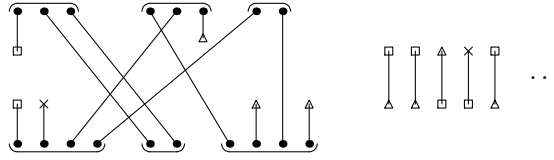}
$$
\caption{Reference  Subsection \ref{ss:young-category}.
 \label{fig:ness-1}}
\end{figure}

 {\bf \punct The category of double cosets.%
 \label{ss:young-category}} Denote
$$K[\alpha_1,\dots,\alpha_m]=K_1[\alpha_1]\times\dots\times K_m[\alpha_m]
.
$$
Consider the category of double cosets
$$
K[\alpha_1,\dots,\alpha_m]\setminus G/K[\beta_1,\dots,\beta_m]
.
$$ 
To be definite set
 $m=3$, see Fig. \ref{fig:ness-1}.
Denote
$$
J_\nu[\gamma]:=\{0,1,\dots,\gamma\}\times \{\nu\}\subset \N\times\{1,2,3\}
.
$$
We draw the diagram with arcs. To the upper row we put
$J_1[\beta_1]\cup J_2[\beta_2]\cup J_3[\beta_3]$, to the lower row 
elements of
$J_1[\alpha_1]\cup J_2[\alpha_2]\cup J_3[\alpha_3]$. We mark different colours ($\nu=1,2,3$) by
roods, squares, and triangles.
Let
 $g$ send  
$(l,\nu)$ to $(k,\mu)$. We draw the following objects depending on
positions of the pairs
 $(l,\nu)$ and  $(k,\mu)$:

$\bullet$ if $(l,\nu)\in J_\nu[\beta_\nu]$,  $(k,\mu)\in J_\nu[\alpha_\mu]$,
then we draw an arc from 
 $(l,\nu)$ to $(k,\mu)$;

$\bullet$ if  $(l,\nu)\in J_\nu[\beta_\nu]$,  $(k,\mu)\notin J_\nu[\alpha_\mu]$,
then we draw an arc with the origin at
 $(l,\nu)$ and a free lower end marked by the color 
 $\mu$; in the same way we proceed in the case
$(l,\nu)\notin J_\nu[\beta_\nu]$,  $(k,\mu)\in J_\nu[\alpha_\mu]$ 
(we draw an arc with the end in  $(l,\mu)$ and a free upper end marked by 
the color 
 $\nu$);
 
$\bullet$ if   $(l,\nu)\notin J_\nu[\beta_\nu]$,  $(k,\mu)\notin J_\nu[\alpha_\mu]$ and
$\mu\ne\nu$, then we draw an arc whose upper end is marked by the colour
$\nu$, and the lower end is marked by the colour  $\mu$
 (we distinguish the upper and lower ends of the arc);
 
$\bullet$ if $(l,\nu)\notin J_\nu[\beta_\nu]$,  $(k,\mu)\notin J_\nu[\alpha_\mu]$ and $\mu=\nu$,
then we draw nothing  (as a result the final picture is finite).
 
 The product is a gluing of diagrams as in the case of chips.

 \sm
 
 {\bf\punct The semigroup   $\Gamma:=K\setminus G/K$ and its characters.%
 \label{ss:char}}
 According the construction of the previous subsection, to a double coset 
 $\in K\setminus G/K$ we assign a collection of arcs  from top to downwards 
 with colored ends. It is important only the number
 of arcs 
 $s_{\nu\mu}$ with each coloring. 
Thus, for any double coset we assign
a collection of numbers
 $s_{\nu\mu}$, where $\nu$, $\mu\le m$ and $\mu\ne\nu$ (i.e., these
 numbers constitute a matrix  $S$ of size   $m\times m$ 
 without the diagonal%
 \footnote{It is possible to set $\infty$ to the diagonal entries.}).
  The numbers $s_{\nu\mu}$ are integer non-negative and satisfy the conditions
  (\ref{eq:s-mu-nu}).
  The multiplication of double cossets corresponds to the addition of matrices.

 As in the case of trisymmetric group, the semigroup
 $\Gamma:=K\setminus G/K]$ has a canonical embedding to each semigroup
 $
 K[\alpha_1,\dots,\alpha_m]\setminus G/K[\alpha_1,\dots,\alpha_m]
 $.
 For each unitary
 $K$-spherical representation of the group  $G$, we have
 a homomorphism  $\chi$ from $K\setminus G/K$ to the multiplicative semigroup
 of complex numbers whose absolute value $\le 1$, 
 $$
 \chi(S_1+S_2)=\chi(S_1)\chi(S_2)
 .$$ 
 Also,
 $$
 \chi(S^t)=\ov {\chi(S)}.
 $$
 
 On our language, Theorem
 \ref{th:ness} is equivalent to the following statement:

 \sm

--- {\it For any unitary  $K$-spherical representation
 $\rho$ of the group  $G$ there is an Hermitian positive semi-definite matrix
 $A$ of order $m$ with units on the diagonal
 such that the spherical character of
 $\rho$ equals 
$$
\chi(S)=\prod_{\mu,\nu:\,\mu\ne\nu} a_{\nu\mu}^{s_{\nu\mu}}
.
$$
} 

This statement would be entirely simple if to admit that the character 
 $\chi$ does not vanish. However, we must examine the structure of the semigroup of its zeros
 and this enlarges a proof.

 Denote by 
 $E_{\nu\mu}$ the matrix-unit, i.e., the matrix 
 having 1 on the place 
 $\nu\mu$, all other matrix elements are zero.
 We will represent matrices $S$ 
 as  $S=\sum s_{\nu\mu} E_{\nu\mu}$. Let as call by {\it cycles} elements of  $\Gamma$ 
 having the form
 $$
 \sigma[l_1\dots l_p]=E_{l_1 l_2} + E_{l_2 l_3}+\dots + E_{l_p l_1},
 \quad \text{where $l_i\ne l_{j}$}
 .$$

\begin{lemma}
The semigroup
 $\Gamma$ is generated by cycles.
\end{lemma}

{\sc Proof.}  Consider a matrix  $S\in \Gamma$. Let $s_{p_1p_2}>0$.
By the condition 
(\ref{eq:s-mu-nu}) there is   $p_3$ such that   $s_{p_2p_3}>0$. If
$p_3\ne p_1$, we take 
$s_{p_3p_4}>0$ etc., until we meet   $p_i=p_j$. Then the matrix 
$S-\sigma[p_i p_{i+1}\dots p_j]$ is non-negative and satisfies 
(\ref{eq:s-mu-nu}).  
\hfill $\square$

\sm

Fix a spherical representation
 $\rho$ and the corresponding character  $\chi$.
We say that colours $\mu$ and $\nu$ {\it are contained in the same component}
if there exists a cycle 
$\sigma$ containing $\mu$, $\nu$ such that  $\chi(\sigma)\ne 0$.

\begin{lemma}
If $\mu$ and $\nu$ are not contained in one component, then
for any matrix $S\in\Gamma$ with $s_{\mu\nu}>0$ {\rm(}or $s_{\nu\mu}>0${\rm)},
we have $\chi(S)=0$.
\end{lemma}
 
{\sc Proof.} It suffices to decompose   $S$ as a sum of cycles.
\hfill $\square$ 

\begin{lemma}
Let $\nu$, $\mu$ be contain in one component.
Then there is a chain 
$\nu_1=\nu$, $\nu_2$, \dots, $\nu_k=\mu$ such that for all $j$
$$
\chi(E_{\nu_j\nu_{j+1}}+E_{\nu_{j+1}\nu_{j}})\ne 0.
$$
\end{lemma}

Note that
$$
E_{\mu\nu}+E_{\nu\mu}=\sigma[\mu\nu].
$$

{\sc Proof.} We choose a cycle   $\sigma[l_1\dots l_p]$, 
containing $\nu$ and $\mu$.
Then
$$\sigma[l_1\dots l_p]+\sigma[l_1\dots l_p]^t=
\sigma[l_1l_2]+\sigma[l_2l_3]+\dots +\sigma[l_pl_1].
$$
Since  $\chi(\sigma+\sigma^t)=|\chi(\sigma)|^2\ne0$, we have  
$\chi(\sigma[l_kl_{k+1}])\ne 0$.
\hfill $\square$ 

\begin{lemma}
Let $\chi(\sigma[l_1l_2])=0$. Then for any cycle 
$\sigma[l_1 l_2\dots l_p]$, we have  $\chi\bigl(\sigma [l_1 l_2\dots l_p]\bigr)=0$.
\end{lemma}

 {\sc Proof.} The matrix
 $(\sigma[l_1 l_2\dots l_p]+\sigma[l_1 l_2\dots l_p]^t)-\sigma[l_1l_2]$
 is contained in  $\Gamma$.
 \hfill $\square$
 
\sm

{\bf \punct Proof of theorem  \ref{th:ness}.%
\label{ss:proof}}
Preserve the notation of the previous subsection. To be definite,
consider the component containing the colour 
 $\nu=1$. Without loss of generality, we can assume that the component
 consists of consecutive colours
 $\nu=1$, 2,\dots, $n$. Consider the graph  $\Delta$,
 whose vertices are the points 
 $1$, 2,\dots, $n$, an edge connects two vertices
 $\mu$,  $\nu$ if 
 $\chi\bigl(\sigma[\mu\nu]\bigr)\ne 0$. Consider an arbitrary spanning tree $\Xi$ of this graph.

\begin{figure}
 $$\epsfbox{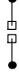}$$
 \caption{Elements $\pi_j$.
 \label{fig:pi}}
 \end{figure}

 Consider the group $K[1,0,\dots,0]$, the subspace   $V$
 of 
 $K[1,0,\dots,0]$-fixed vectors and the representation of the semigroup 
 $$\Gamma_1:=K[1,0,\dots,0]\setminus G/ K[1,0,\dots,0]$$
  in $V$. Consider the  elements  $\pi_j$ of the semigroup 
  determined by pictures shown  on Fig.
 \ref{fig:pi}.
 It is easy to see that 
 $\pi_j^2=\pi_j$. Therefore the corresponding operators $P_j$
 are orthogonal projectors   in $V$.

Notice that 
  $$P_i P_j P_i= P_i\cdot \chi\bigl(\sigma[ij]\bigr).$$
In the left-hand side we have a positive operator, therefore 
  $$
  \chi(\sigma(ij))\ge 0.
  $$
  
  Next, we write a more general identity. For a cycle 
 $\sigma[\nu_1\nu_2\dots\nu_q]$ we have 
\begin{equation}
P_{\nu_q} P_{\nu_1} P_{\nu_2}\dots P_{\nu_q}= \chi\bigl(\sigma[\nu_1\nu_2\dots\nu_q]\bigr)
P_{\nu_q}
\label{eq:proizvedenie}
\end{equation}
 
\begin{lemma}
For each  $j\le n$ the image of the projector  $P_j$ is one-dimensional.
\end{lemma}

 {\sc Proof.} The operator  $P_1$ is the projector to the spherical vector
 (in the notation of Subsection \ref{ss:extension-semigroup}
 it is the matrix  0 in the semigroup of  1-0-matrices
   of size $1\times 1$).  
   Let $\chi\bigl(\sigma[1j]\bigr)=s\ne 0$. Then $P_1 P_j P_1=s P_1$, therefore
   $P_j\ne 0$. On the other hand, $P_j P_1 P_j=s P_j$, 
   therefore rank of  $P_j$ does not exceed  1.
 Thus, for all   $j$ adjacent to 1 in the graph  $\Delta$, we have  $\rk P_j=1$.
 Next, we repeat the same reasoning  for $P_j$  instead of $P_1$
 and refer to the connectedness of the graph 
 $\Delta$.
    \hfill $\square$

\sm

Denote by 
 $L_\nu$ the line   $P_\nu V$. Choose in  $L_1$ a unit vector 
$v_1$. Let $\nu_1=1$, $\nu_2$, \dots, $\nu_q=\nu$ be a way from  1 to $\nu$ along the 
spanning tree
 $\Xi$.
Consider the vector
$$
w_\nu:=P_{\nu_q}\dots P_{\nu_1}v_1\in L_\nu
.
$$
This vector is non-zero, moreover,
\begin{multline*}
\la w_\nu,w_\nu\ra=\la P_{\nu_q}\dots P_{\nu_1}v_1,P_{\nu_q}\dots P_{\nu_1}v_1\ra=
\la P_{\nu_1}\dots P_{\nu_q}\dots P_{\nu_1}v_1,v_1\ra
=\\=
\la v_1,v_1\ra \prod\nolimits_q \chi\bigl(\sigma[\nu_q \nu_{q+1}]\bigr)
.
\end{multline*}
Set
$$
v_\nu:={w_\nu}/{\|w_\nu\|}.
$$
Then for all pairs
 $\mu$, $\nu$, which are adjacent in the spanning tree,
we have
\begin{equation}
\la v_\mu, v_\nu\ra=\chi\bigl(\sigma[\mu\nu]\bigr)^{1/2}
.
\label{eq:-1/2}
\end{equation}
Let $S$ be the matrix of size  $n\times n$ with non-negative integer elements
(as above, diagonal elements are not defined,
we do not impose  conditions   (\ref{eq:s-mu-nu})).
Set
$$
\psi(S)=\prod\nolimits_{\mu\ne \nu} \la v_\mu, v_\nu\ra^{s_{\mu\nu}}
.
$$

\begin{lemma}
If  $S$ satisfies the conditions
\begin{equation}
\sum\nolimits_{\mu: \, \mu\le n, \mu\ne \nu} s_{\mu\nu}= 
\sum\nolimits_{\mu: \, \mu\le n, \mu\ne \nu} s_{\nu\mu} 
,
\label{eq:psi-chi}
\end{equation}
 then $\psi(S)=\chi(S)$.
\end{lemma}
   
 {\sc Proof.} Recall that the semigroup $\Gamma_{(1)}$ of all  $S$,
 satisfying
 (\ref{eq:psi-chi}) is generated by cycles.
If a cycle does not contained in the graph, then in the both sides 
of 
(\ref{eq:psi-chi}) we have zero. 
If a cycle has the form
 $\sigma[\mu\nu]$,
where   $\mu\nu$ is an edge of the spanning tree, then
both sides coincide by  (\ref{eq:-1/2}).
 Let $\mu\nu$ be an edge of the graph  $\Delta$, which does not contained in the tree
 $\Xi$.
Adding  $\mu\nu$ to $\Xi$ we get the graph with a unique cycle.
Pass the cycle in the direction
 $\mu\nu$ and take the corresponding element of 
 $\sigma[\kappa_1\dots \kappa_p]$
 where $\kappa_1=\nu$, $\kappa_p=\mu$ of the semigroup  
 $\Gamma_{(1)}$. For this case the identity  (\ref{eq:psi-chi}) 
 is reduced to (\ref{eq:proizvedenie}).
 
 Now, consider an arbitrary cycle 
 $C$ passing along edges of $\Delta$.
 It as a linear combination
 with integer elements of cycles of the form 
 $\sigma[\mu\nu]$, 
 $ \sigma[\kappa_1\dots \kappa_p]$ considered above.
Keeping in the mind that values of  $\psi$ on the cycles $\sigma[\mu\nu]$, 
 $ \sigma[\kappa_1\dots \kappa_p]$
 are invertible, we get that the statement is valid for 
 $C$.
   \hfill $\square$
 
 \sm
 
 Applying this lemma to all components of the set
 $\{1,\dots,m\}$ we get Theorem 
 \ref{th:ness}.    \hfill $\square$

 \sm
 
 {\bf \punct Remarks.} As we have seen above,
 there are many spherical pairs  $G\supset K$,
related to infinite symmetric groups.
Completeness of list of spherical functions 
(see conjectures from  \ref{ss:conjectures}) 
now is proved only in the following cases,
the Thoma theorem, the Nessonov theorem, and also in 3 cases 
relative to bisymmetric group, they were considered by Olshanski in
 \cite{Olsh-symm}: 
 
 \sm
 
 1. $G=S_{2\infty}$ and $K$ is the hyperoctahedral group  $S_\infty \wre \Z_2$.
 
 \sm
 
 2. $G=S_{2\infty+1}$, $K=S_\infty \wre \Z_2$.
 
 \sm
 
 3. $G=S_{\infty+1}\times S_\infty$, $K=\diag S_{\infty}$.

 Passing to the language of train categories, we observe that these 
 pairs are really simpler than other pairs. The corresponding combinatorial
 objects are one-dimensional, and additional structure attributed to arcs 
 is a number (a 'length' in  \cite{Olsh-symm}) and
 not a color tiling as in general constructions.

There are other cases possesing this property. One of them 
is the symmetric group
 $S_{3\infty}$ with the subgroup  $ (S_\infty \wre \Z_2)\times S_\infty$ 
 (we also can consider
 $S_{4\infty}\supset (S_\infty \wre \Z_2)\times S_\infty\times S_\infty$
 etc., but two hyperoctahedral factors is a complicated case).

Next level of complexity corresponds to spherical pairs associated with chips,
the first is the example from Subsection
\ref{ss:chipy-strelki}
and a pair considered by Nessonov in  \cite{Ness-fa}.

\noindent
\tt Math.Dept., University of Vienna,
 \\
 Oskar-Morgenstern-Platz 1, 1090 Wien;
 \\
\& Institute for Theoretical and Experimental Physics (Moscow);
\\
\& Mech.Math.Dept., Moscow State University.
\\
\& Institute for Information Transmisiion, Moscow
\\
e-mail: neretin(at) mccme.ru
\\
URL:www.mat.univie.ac.at/$\sim$neretin

\end{document}